\documentclass[12pt]{article}
\usepackage{amsmath}
\usepackage{amssymb}

\newtheorem{conjecture}{Conjecture}[section]

\newtheorem{lemma}[conjecture]{Lemma}
\newtheorem{proposition}[conjecture]{Proposition}

\newtheorem{theorem}[conjecture]{Theorem}

\newenvironment{proof}{\noindent {\bf Proof} \hspace{.1cm}}{\hfill
${\bf QED}$ \\ \vspace{.15cm}}

\title{Conformal measures associated to ends of hyperbolic $n$-manifolds}
\author{James W. Anderson \\
Kurt Falk\footnote{Research supported by the V\"ais\"al\"a Foundation
and the University of Helsinki} \\
Pekka Tukia}
\date{\today}

\begin{document}

\maketitle

\abstract{
\noindent
Let
$\Gamma$ be a non-elementary Kleinian group acting on the closed
$n$-dimensional unit ball and assume that its Poincar\'e
series converges at the exponent $\alpha$.  Let $M_\Gamma$ be the
$\Gamma$-quotient of the open unit ball.
We consider certain families ${\cal E}=\{E_1,...,E_p\}$ of open subsets
of $M_\Gamma$  such that $M_\Gamma\setminus(\cup_{E\in {\cal E}} E)$
is compact. The sets  $E_i$ are called ends of $M_\Gamma$ and ${\cal E}$ is
called a complete collection of ends for $M_\Gamma$.
We show that we can associate to each end $E\in{\cal E}$
a conformal measure of dimension $\alpha$ such that the 
measures corresponding to
different ends are mutually singular if non-trivial. Each conformal
measure for $\Gamma$ of dimension $\alpha$ on the limit set $\Lambda(\Gamma)$ of $\Gamma$ 
can be written as a sum of such conformal
measures associated to ends  $E\in{\cal E}$. In dimension 3, our results
overlap with some results of Bishop and Jones \cite{bijo}.
}


\section{Introduction and survey of known results}
\label{introduction}

\medskip
\noindent
Consider a non-elementary Kleinian group $\Gamma$ acting on the closed
$n$-dimensional unit ball ${\bf B}^n\cup {\bf S}^{n-1}$ and assume that its Poincar\'e
series converges at the exponent $\alpha$.  Let $M_\Gamma$ be the
$\Gamma$-quotient of the open unit ball ${\bf B}^n$.
In this note we study certain families ${\cal E}=\{E_1,...,E_p\}$ 
of open subsets of $M_\Gamma$  such that $M_\Gamma\setminus (\cup_{E\in {\cal E}} E)$
is compact. We call the sets  $E_i$ ends of $M_\Gamma$, and ${\cal E}$ 
a complete collection of ends for $M_\Gamma$.
A point $z$ on the unit sphere ${\bf S}^{n-1}$ is called an endpoint for an end $E$ if
any geodesic ray $R$ in ${\bf B}^n$ towards $z$ contains a subray which projects 
onto a ray $R'$ in $M_\Gamma$ so that $R' \subset E$ and so that
the distance of a point $x$ on $R'$ to the boundary
of $E$ tends to infinity as $x$ tends to $z$ along $R'$. The point $z$ 
is called an end limit point if $z$ is a limit point of $\Gamma$. 
We construct and investigate conformal measures concentrated on the set 
of endpoints for an end. 

\medskip
\noindent
We sometimes refer to a {\em conformal measure of dimension $\alpha$ for $\Gamma$} as an {\em $\alpha$-conformal measure for $\Gamma$}.  Note that we do not assume that a 
conformal measure for $\Gamma$ is necessarily supported by the limit set 
$\Lambda(\Gamma)$.  The main result of this paper, stated and proved in the text as Theorem \ref{main-theorem}, can be paraphrased as follows.

\medskip
\noindent
{\bf Theorem}
{\em
Let $\Gamma$ be a non-elementary Kleinian group $\Gamma$ acting on ${\bf B}^n$
which has a complete collection of ends ${\cal E}=\{E_1,\ldots, E_p\}$ 
and assume that its
Poincar\'e series converges at $\alpha$.
Let $\Lambda_i$ be the set of end limit points of $E_i$.  If  
$\Lambda_i\neq\emptyset$, 
there exists a non-trivial 
$\alpha$-conformal measure for $\Gamma$
supported  by the end limit points of $E_i$, 
and any two such measures
corresponding to different ends are mutually singular. Each
$\alpha$-conformal measure $m$ for $\Gamma$ on $\Lambda(\Gamma)$ 
can be written as a sum of
$\alpha$-conformal measures $m_i$ for $\Gamma$ 
supported by the $\Lambda_i$.  
}

\medskip
\noindent
As a general result, it is known, see Sullivan \cite{sullivan1}, Roblin
\cite{roblin}, that for a Kleinian group $\Gamma$ acting on ${\bf B}^n$
of divergence type, there exists a unique invariant conformal measure on
its limit set, up to multiplication by constants.  Specializing to the
case of $n=3$, suppose that $\Gamma$ is a topologically tame Kleinian
group acting on ${\bf B}^3$ for which $\Lambda(\Gamma) ={\bf S}^2$.
Then, every $\Gamma$-invariant conformal measure on ${\bf S}^2$ is a
multiple of Lebesgue measure on ${\bf S}^2$.  In this generality, this
result can be obtained by combining Theorem 9.1 of Canary \cite{canary2}
with Proposition 3.9 of Culler and Shalen \cite{culler-shalen}.  Note that
while Proposition 3.9 of \cite{culler-shalen} holds for all ${\bf B}^n$,
Theorem 9.1 of \cite{canary2} is specific to ${\bf B}^3$.

\medskip
\noindent
Therefore we need to consider only groups of convergence type.  It is known
that a topologically tame Kleinian group $\Gamma$ acting on ${\bf B}^3$ is
of convergence type if and only if $\Lambda(\Gamma)\ne {\bf S}^2$, see
Corollary 9.9.3 of Thurston \cite{thurston}. We note that there are
Kleinian groups of convergence type acting on ${\bf B}^3$ for which the
invariant conformal measure on its limit set is unique up to
multiplication by constants, such as the examples given by Sullivan
\cite{sullivan-growth}, which inspired the present work.

\medskip
\noindent
Bishop and Jones, see Corollary 1.3 of \cite{bijo}, prove that if $\Gamma$ is
a topologically tame, geometrically infinite Kleinian group acting on
${\bf B}^3$ for which the injectivity radius of
$M_\Gamma ={\bf B}^3/\Gamma$ is bounded away from zero and
$\Lambda(\Gamma) \ne {\bf S}^2$, then ${\cal H}^\varphi$,
the Hausdorff measure associated to the gauge function
$\varphi(t) =t^2\sqrt{\log(\frac{1}{t}) \log\log\log(\frac{1}{t})}$,
is a conformal density of dimension $2$.
As a consequence they obtain that for each
geometrically infinite end there is a unique $2$-conformal measure
(up to multiplicative constants), the action of $\Gamma$ on the
boundary is ergodic with respect to each of these measures, these
measures are mutually singular, and finally, any $2$-conformal measure for
$\Gamma$ is a linear combination of them.
The methods they use are analytic, constructing a positive harmonic
function on $M_\Gamma$ which grows at most linearly in the geometrically
infinite ends of $M_\Gamma$.

\medskip
\noindent
If applied to the situation in $3$ dimensions, Theorem \ref{main-theorem} 
means that one can remove the lower bound on the injectivity radius of
$M_\Gamma$ and replace the condition that the group is of the 
second kind by the convergence of the Poincar\'e series
in Bishop and Jones' result.  However, our statements are not as strong as theirs, in that we are unable to prove ergodicity.  Our arguments employ
straightforward topological and dynamical mechanisms, and we feel that the straightforward nature of the proofs makes up for this lack of ergodicity.  It is also not clear the extent to which one can expect such strong ergodicity results to hold in all dimensions.

\medskip
\noindent
A note on referencing: it is our intention to give due credit to all
authors, but we have chosen to sometimes reference standard texts rather
than the first statement of a result.  As a standard reference on the
basics of Kleinian groups, we use Maskit \cite{maskit}.  As a standard
reference on the measure theoretic constructions involving Kleinian
groups, we use Nicholls \cite{nicholls}.

\medskip
\noindent
This paper is structured as follows. In Section~\ref{basic-definitions} we
give the basic definitions and fix the notation. 
In Section~\ref{ends-and-endgroups} we introduce
the notions of an end of a hyperbolic $n$-manifold and the associated end
groups as they are used in this note, and give the topological and
geometric properties which are necessary for the measure-theoretical
considerations of Section~\ref{measures}. 
In Section \ref{measures}, we state and prove our main
results and discuss some applications.   Finally, in
Section~\ref{ends-and-endgroups-3dim} we specialize the discussion to the
case of a Kleinian group acting on ${\bf B}^3$.

\section{Basic definitions}
\label{basic-definitions}

\medskip
\noindent
Throughout this note, we work in the {\em Poincar\'e ball model} of hyperbolic
$n$-space, where $n\ge 2$.  The underlying space is the unit ball
${\bf B}^n =\{ x \in {\bf R}^n\: |\: |x| <1 \}$ in ${\bf R}^n$,
with the element of arc-length $\frac{2}{1-{|x|}^2} |{\rm d}x|$.  The
hyperbolic distance between points $x$, $y$ in ${\bf B}^n$ is denoted
${\rm d}(x,y)$.  For subsets $X$ and $Y$ of ${\bf B}^n$,
set ${\rm d}(X,Y) =\inf\{ {\rm d}(x,y) \: |\: x\in X, y\in Y\}$.
The {\em sphere at infinity} of ${\bf B}^n$ is the unit
sphere ${\bf S}^{n-1}$ in ${\bf R}^n$.

\medskip
\noindent
 We denote the closure of $X$ in the Euclidean topology on
the closed $n$-ball ${\bf B}^n\cup {\bf S}^{n-1}$ by $\overline{X}$; in
particular, $\overline{{\bf B}^n} = {\bf B}^n\cup {\bf S}^{n-1}$.  
We denote the boundary of $X\subset {\bf B}^n$ in ${\bf B}^n$ by
$\partial X$.
For a
subset $X$ of $\overline{{\bf B}^n}$, let
$\partial_\infty (X) =\overline{X}\cap {\bf S}^{n-1}$.

\medskip
\noindent
A {\em Kleinian group} is a discrete subgroup $\Gamma$ of the group of (possibly orientation-reversing) isometries of hyperbolic $n$-space ${\bf B}^n$.  (We note that Kleinian groups are often assumed to contain only orientation-preserving isometries.  This is the case in many of the papers we have referred to.  However, this assumption is not relevant to our arguments in Sections \ref{ends-and-endgroups} and \ref{measures}, and so we do not make it here.)  A Kleinian group is {\em elementary} if it contains an abelian
subgroup of finite index, and is {\em non-elementary} otherwise.  Unless otherwise stated, we will assume that a Kleinian group is non-elementary.  We denote the induced hyperbolic distance between points $x$ and $y$ in the quotient ${\bf B}^n/\Gamma$ by ${\rm d}(x,y)$.

\medskip
\noindent
An orientation-preserving isometry of ${\bf B}^n$ extends to a conformal
homeomorphism of the sphere at infinity ${\bf S}^{n-1}$.  The {\em domain
of discontinuity} $\Omega(\Gamma)$ of a Kleinian group $\Gamma$ is the
largest open subset of ${\bf S}^{n-1}$ on which $\Gamma$ acts properly
discontinuously.  The complement of $\Omega(\Gamma)$ in ${\bf S}^{n-1}$ is
the {\em limit set} $\Lambda(\Gamma)$ of $\Gamma$.  We equivalently define
$\Lambda(\Gamma)$ to be the set of accumulation points of the orbit
$\Gamma x$ for any point $x \in \overline{{\bf B}^n}$.  A Kleinian group
$\Gamma$ is {\em of the first kind} if $\Lambda(\Gamma) ={\bf S}^{n-1}$,
and is {\em of the second kind} otherwise.  In the latter case,
$\Lambda(\Gamma)$ is a closed, nowhere dense subset of ${\bf S}^{n-1}$.
If $\Gamma$ is non-elementary, then $\Lambda(\Gamma)$ is perfect.

\medskip
\noindent
Let $\Gamma$ be a Kleinian group acting on ${\bf B}^n$.  
For a subset $X$ of $\overline{{\bf B}^n}$, the {\em stabilizer} 
$\Gamma_X$
of $X$ in $\Gamma$ is defined to be the subgroup 
\[ \Gamma_X  := \{ \gamma\in\Gamma\: |\: \gamma(X) =X\} \]
of $\Gamma$.  

\medskip
\noindent
A limit point $x\in \Lambda(\Gamma)$ of a Kleinian group $\Gamma$ is a
{\em conical limit point} if the following holds: there exists a
hyperbolic ray $R$ in ${\bf B}^n$ ending at $x$, a number
$\varepsilon>0$, a point $z\in {\bf B}^n$, and a sequence
$\{ \gamma_n\}$ of distinct elements of $\Gamma$ so that
$\gamma_n (z)\in U_\varepsilon(R)$ for all $n$, and
$\gamma_n (z)\rightarrow x$ in $\overline{{\bf B}^n}$. Here,
$U_\varepsilon (R) =\{ z\in {\bf B}^n\: |\: {\rm d}(x,R)<\varepsilon\}$
is the open $\varepsilon$-neighborhood of $R$ in ${\bf B}^n$.
Equivalently, $x$ is a conical limit point of $\Gamma$ if there
exists a sequence $\{ x_n\}$ of points of $R$ converging to $x$ in
$\overline{{\bf B}^n}$ so that, if
$\pi: {\bf B}^n\rightarrow {\bf B}^n/\Gamma$
is the covering map, then the $\pi(x_n)$ all lie in a compact subset of
${\bf B}^n/\Gamma$.  The collection of all conical limit points of
$\Gamma$ is denoted $\Lambda_c(\Gamma)$.

\medskip
\noindent
Let $\Gamma$ be a Kleinian group acting on ${\bf B}^n$.
Fix a point $y\in {\bf B}^n$.  For $x \in {\bf B}^n$, define the {\em
Poincar\'e series} $P_\Gamma(x,y,s)$ of $\Gamma$ based at $x$ to be
\[
P_\Gamma(x,y,s) = \sum_{\gamma\in \Gamma} \exp(-s\, {\rm d} (x,
\gamma(y))).
\]
By the triangle inequality, if $P_\Gamma(x,y,s)$ converges at $s$ for some
$x\in {\bf B}^n$, then it converges at $s$ for all $x\in {\bf B}^n$. The 
same holds for  $y$. 

\medskip
\noindent
The {\em critical exponent} $\delta(\Gamma)$ of $\Gamma$ is defined to
be
\[ \delta = \delta(\Gamma) := 
\inf \{ s>0\: |\: P_\Gamma(x,y,s)\mbox{ converges}
\}. \]
$\Gamma$ is {\em of $\delta$-convergence type}, or simply
{\em of convergence type}, if $P_\Gamma(x,y,\delta)$ converges, and
{\em of divergence type} if $P_\Gamma(x,y,\delta)$ diverges.

\medskip
\noindent
Let $\Gamma$ be a Kleinian group acting on ${\bf B}^n$.
Then, we have that $\delta(\Gamma) \leq n-1$, see Nicholls \cite{nicholls},
Theorem 1.6.1.  If  $\Gamma$ is of the second kind, then $\Gamma$ is of
$(n-1)$-convergence type, see Nicholls \cite{nicholls}, Theorem 1.6.2.
If $\Gamma$ has finite volume quotient ${\bf B}^n/\Gamma$,
then $\delta(\Gamma) = n-1$ and $\Gamma$ is of $(n-1)$-divergence type.
If $\Gamma$ is geometrically finite and of the second kind, then $\delta(\Gamma) < n-1$,
see Sullivan \cite{sullivan2}, Tukia \cite{tukia2}.
We have that $\delta(\Gamma) > 0$ for any non-elementary Kleinian 
group $\Gamma$, see Beardon \cite{beardon}.  

\medskip
\noindent
There exist finitely generated, geometrically infinite
Kleinian groups of the first kind acting on ${\bf B}^3$ which are of $2$-divergence type
(see e.g. Thurston \cite{thurston}, Sullivan \cite{sullivan4},
Rees \cite{rees1}, \cite{rees2},
or Aaronson and Sullivan \cite{aaronson-sullivan}).

\medskip
\noindent
Throughout we shall work with real-valued, non-negative, finite measures
on $\overline{{\bf B}^n}$. We shall say that a {\em support} of such a 
measure is a measurable, not necessarily uniquely determined subset of 
$\overline{{\bf B}^n}$ whose complement has measure zero.

\section{Ends and endgroups of Kleinian groups acting on ${\bf B}^n$}
\label{ends-and-endgroups}

\medskip
\noindent
Let $\Gamma$ be a Kleinian group acting on ${\bf B}^n$.
We associate to $\Gamma$ the following orbit spaces:
\begin{itemize}
\item[] $M_\Gamma ={\bf B}^n/\Gamma$,
\item[]  $\overline{M_\Gamma} =({\bf B}^n\cup \Omega(\Gamma))/\Gamma$,
\item[] $\partial_\infty M_\Gamma =\Omega(\Gamma)/\Gamma$.
\end{itemize}

\medskip
\noindent
For a subset $X$ of $\overline{M_\Gamma}$, 
let $\overline{X}$ denote
the closure of $X$ in the Euclidean topology on $\overline{M_\Gamma}$, and
let $\partial_\infty (X) = \overline{X}\cap \partial_\infty(M_\Gamma)$.

\medskip
\noindent
We now give the basic definition of this paper.  A connected open subset
$E$ of $M_\Gamma$ is an {\it end} of $M_\Gamma$ if $\partial E$ is compact and non-empty 
and if $E$ has non-compact closure in $M_\Gamma$.  
The bulk of the paper is devoted to exploring this definition of end, and
to constructing measures associated to ends. 
Note that the word end has been used in several different senses for 
hyperbolic manifolds.  For instance, the end in Bonahon \cite{bonahon} 
can be described as a certain type of restricted equivalence class of ends in our sense.

\medskip
\noindent
One of the cases we consider in detail, in Section \ref{ends-and-endgroups-3dim},  is that $n=3$ and $\partial E$ is a
separating compact surface, which we will often take to be a boundary
component of a compact core of $M_\Gamma$.  While it
is sometimes the case that $\partial E$ is {\em incompressible}, meaning
that $\pi_1(\partial E)$ is infinite and the inclusion of $\partial E$
into $M_\Gamma$ induces an injection on fundamental groups, the more complicated
case occurs when $\partial E$ is not incompressible.  

\medskip
\noindent
Let $\pi :{\bf B}^n\to M_\Gamma$ be the covering projection, let $E$ be an
end of $M_\Gamma$, and let $E^0_i$, $i\in I$, be the components of
$\pi^{-1}(E)$. Let 
$\Gamma_i:= \Gamma_{E^0_i}$ 
be the stabilizer
of $E^0_i$ in $\Gamma$. We call $E^0_i$ an {\it end} of $\Gamma$ and
$\Gamma _i$ the corresponding {\it end group} of $\Gamma$.  
By construction,
the groups $\Gamma _i$ are conjugate subgroups of $\Gamma$. 

\medskip
\noindent
If we say that $\widetilde E$ is an end of $\Gamma$, we mean that  
$\widetilde E$ is obtained 
as above, i.e. $\widetilde E$ is a component of $\pi^{-1}(E)$ for an end  
$E$ of $M_\Gamma$. 
Thus ends of $\Gamma$ are subsets of ${\bf B}^n$, 
while the ends of $M_\Gamma$ are subsets of $M_\Gamma$. 

\medskip
\noindent
Let
$F^0_i = 
\big( \, \overline{{\bf B}^n} \setminus \overline{E^0_i} \, \big)
\cup \partial E^0_i$.
Since the sets $F^0_i/\Gamma _i$ are canonically homeomorphic, we set
$F=F^0_i/\Gamma_i$; in some loose sense, $F$ captures the behavior of an
end group on the complement of its corresponding end.  Note that $F^0_i$
is invariant under the action of $\Gamma_i$, since
$\overline{{\bf B}^n}$, $\overline{E^0_i}$, and $\partial E^0_i$ are all
invariant under the action of $\Gamma_i$.  Also, since 
$\overline{E^0_i}\cap {\bf S}^{n-1}$ is a non-empty closed subset of 
${\bf S}^{n-1}$ invariant under $\Gamma_i$, we see that
$\Lambda(\Gamma_i)\subset \overline{E^0_i}\cap {\bf S}^{n-1}$, and so
$F^0_i\subset {\bf B}^n\cup \Omega(\Gamma_i)$.  If $F$ is compact,
we say that $E$ is a {\it bounded end}.  Additionally, in this case we
refer to the $E^0_i$ as {\it bounded ends} of $\Gamma$ and the $\Gamma _i$
as {\it bounded end groups}. Thus, if we refer to $\widetilde E$ as a 
bounded end of $\Gamma$, $\widetilde E$ is obtained in this manner.

\medskip
\noindent
We can equivalently characterize a bounded end as follows: Since each
$E^0_i$ is precisely invariant under its stabilizer $\Gamma_i$ in
$\Gamma$, we can identify the quotient $E^0_i/\Gamma _i$ with the end $E$
of $M_\Gamma$, and so we can regard $E$ as a subset of
$\overline{M_i}=\overline{M_{\Gamma _i}}$. The end $E$ of $M_\Gamma$
(or equivalently the end $E^0_i$ of $\Gamma$) is a bounded end if
$\overline{M_i} \setminus (E\cup\partial_\infty E)$ is compact.

\medskip
\noindent
Let $\Gamma$ be a Kleinian group acting on ${\bf B}^n$.  We say that a
point $z\in {\bf S}^{n-1}$ is an {\it endpoint} of an end $E^0$ of
$\Gamma$ if, whenever $R$ is a hyperbolic ray with endpoint at infinity
$z$, there exists a subray $R'$ of $R$ that is contained in $E^0$, for
which the hyperbolic distance ${\rm d}(x,\partial E^0)\to \infty $ as
$x\to z$ on $R'$. Note that if  $S'$  is the projection of  $R'$ to 
the quotient $M_\Gamma$, then we still have that ${\rm d}(x,\partial E)\to\infty$ as  $x$ 
tends towards infinity on  $S'$ and where $E^0$ projects to  $E$.
Note that by definition, an endpoint of an
end $E^0$ of $\Gamma$ is never a conical limit point of $\Gamma$.

\medskip
\noindent
If $\Phi$ is the end group associated to $E^0$, that is
$\Phi = \Gamma_{E^0}$,
if $z$ is an endpoint of $E^0$, and if
in addition $z\in \Lambda(\Phi)$, we say that $z$ is an
{\it end limit point} of $E^0$. We will see (Lemma \ref{tukia-lemma-1})
that $z\in \Lambda(\Phi)$ as 
soon as $z\in\Lambda(\Gamma)$. We denote the set of end limit points
of $\Phi$ and $E^0$ by $\Lambda_e(E^0) = \Lambda_e(E^0, \Phi)$.
Note that we have the inclusion
$\Lambda_e(E^0, \Phi)\subset \Lambda(\Phi) \setminus \Lambda_c(\Phi)$.  
If the end is bounded, this is an equality (cf. Lemma 
\ref{tukia-lemma-1}) but this need not be so in general.

\medskip
\noindent
{\sc Example 1.} 
To give a concrete example of an end, consider the following example
in the case $n=2$. Let $S$ be a Riemann surface of genus $0$ and infinite
analytic type; for instance, let $S$ be the domain of discontinuity of a
$2$-generator Schottky group $\Phi$.  We can uniformize $S$ by a Fuchsian
group $\Gamma$, so that $S ={\bf B}^2/\Gamma$.  Let $c$ be a simple closed
geodesic on $\Omega(\Phi)/\Phi$ which lifts to a simple closed geodesic
$C$ on $\Omega(\Phi) =S$ (such a curve always exists), and consider the
curves $\{ \varphi(C)\: | \: \varphi\in \Phi\}$ on $S$.  Each of these
curves is compact and separating, and so each determines a pair of ends,
namely the two components of $S \setminus \varphi(C)$.  Now,  let $C_1$ and $C_2$ be two disjoint
lifts of the curve $c$ on $\Omega(\Phi)/\Phi$ to $S =\Omega(\Phi)$.  For
$C_1$, we can choose the end $E_1$ containing $C_2$, and for $C_2$, we can
choose the end $E_2$ contained in $E_1$.  Obviously, we can continue 
this and obtain an infinite sequence of ends  
$E_1\supset E_2\supset...$ and there is no smallest end in the sequence.  
The nested sequence of 
ends induces an inclusion on the sets of end limit points, that is, 
$\Lambda_e(E_1)\supset \Lambda_e(E_2)\supset...\,\,$.

\medskip
\noindent
{\sc Example 2.} 
The following is a simple example of a bounded end.  Let $\Gamma$ be a
Kleinian group acting on ${\bf B}^n$.  Let $v$ be a
parabolic fixed point of $\Gamma$ such that 
$\Gamma_{v}$
is free abelian of rank $n-1$. There then
exists an open horoball $B$ at $v$ that is precisely invariant under $\Gamma_{v}$
in $\Gamma$, so that $\gamma(B) =B$ for all $\gamma\in \Gamma_{v}$ and
$\gamma(B)\cap B=\emptyset$ for all $\gamma\in \Gamma \setminus \Gamma_{v}$.
(As usual, an open horoball is a Euclidean ball contained in 
${\bf B}^n$ whose boundary sphere is tangential to ${\bf S}^{n-1}$.)
Hence, the stabilizers 
$\Gamma_{B}$ and $\Gamma_{v}$
of $B$ and $v$, respectively, coincide. Let
$S^0=\partial B$, and note that $\overline{B}$ is the disjoint union of $B$,
$\partial B$, and $\partial_\infty(B) =\{ v\}$. Then
$S=S^0/\Gamma =S^0/\Gamma_{v}$ is compact since $\Gamma_{v}$ has full rank.
The same is true of $(\overline{{\bf B}^n}\setminus (B\cup \{v\}))/\Gamma$,
which is homeomorphic to $S\times [0,1]$. Thus $B$ is a bounded end of
$\Gamma$ with end group $\Gamma_{v}$, and $B/\Gamma$ is a bounded end of
$M_\Gamma$. Let $R$ be any geodesic ray in ${\bf B}^n$ ending at $v$.
Then, it follows from basic properties of hyperbolic space that there
exists a subray $R'$ of $R$ that is contained in $B$, and for which the
hyperbolic distance ${\rm d}(x,\partial B)\to \infty $ as $x\to z$ on
$R'$.  Hence, $v$ is an end limit point of $B$.  In fact, in this case
we have that $\{ v\} =\Lambda_e(B) =\Lambda(\Gamma_{v}) \setminus\Lambda_c(\Gamma_{v})$,
since any hyperbolic ray ending at any point $z\ne v$ in ${\bf S}^{n-1}$
must eventually exit any given horoball based at $z$.

\medskip
\noindent
The crucial fact about bounded ends is the following tripartite division
of points of the sphere at infinity ${\bf S}^{n-1}$ of ${\bf B}^n$.

\begin{lemma}
Let $\Gamma$ be a Kleinian group acting on ${\bf B}^n$.
Let $\Phi$ be a bounded end group of $\Gamma$ associated to the 
bounded end $E^0$ of $\Gamma$. Let $z\in {\bf S}^{n-1}$. 
Then either $z\in \Omega (\Phi)$,
$z\in \Lambda _c(\Phi)$, or $z$ is an end limit point of $E^0$. 

\medskip
\noindent
If $z\in \Lambda (\Gamma)$ and $z$ is an endpoint of $E^0 $, then
$z\in \Lambda (\Phi)$ and this is true even if the end  $E^0$ is not 
bounded.
\label{tukia-lemma-1}
\end{lemma}

\begin{proof}
Let $F^0 = (\overline{{\bf B}^n} \setminus \overline{E^0}) \cup \partial E^0$, 
and let $F=F^0/\Phi$.
Let $D = F^0\cap {\bf S}^{n-1}= {\bf S}^{n-1}\setminus \overline{E^0}$.
By definition of a bounded end, $F^0/\Phi$ is compact, and hence
$\partial_\infty (F^0/\Phi) = D/\Phi$ is compact as well.

\medskip
\noindent
Since $D/\Phi$ is compact, there are finitely many open hyperbolic
half-spaces $H_j$, $1\le j\le p$, contained in ${\bf B}^n\setminus E^0$ so
that $W=\cup_{\varphi \in \Phi, 1\le j\le p} \: \varphi (H_j)$ is a
neighborhood of $D$ in $\overline{{\bf B}^n}$.  The quotient
$(F^0\setminus W)/\Phi $ is compact, as it is a closed subset of the
compact set $F^0/\Phi$.

\medskip
\noindent
Let $z\in {\bf S}^{n-1}$ and let $R$ be a hyperbolic ray with endpoint at
infinity $z$. If $R$ has a subray $R'$ with endpoint at infinity $z$ such
that $R'\subset E^0$, then as we move $x$ along $R'$ towards $z$, there
are two possibilities.  One is that ${\rm d}(x,\partial E^0)\to \infty$ as
$x\rightarrow z$, and hence $z$ is an endpoint of $E^0$; in this case,
either $z\in \Omega(\Phi)$ or $z\in \Lambda _e(E^0)$.  The other
case is that there is a positive number $r$ so that ${\rm d}(x,\partial
E^0)<r$ for $x\in R'$ arbitrarily close to $z$. Since $\partial
E^0/\Gamma$ is compact, it follows that in this latter case
$z\in \Lambda_c(\Phi)$.

\medskip
\noindent
If there is no subray of $R$ contained in $E^0$, then let
$X=F^0\setminus W$ (with $W$ constructed as above) so that $X/\Phi $ is compact.  In this case, either $R$
contains points $x_i\in R\cap X$ such that $x_i\to z$, in which case
$z\in \Lambda _c(\Phi)$, or else $R$ meets some $\gamma H_i$,
$\gamma \in \Phi$. But then $R$ contains a subray $R'$ with endpoint at
infinity $z$ such that $R'\subset \gamma H_i$ and it follows that
$z\in \Omega(\Phi )$.

\medskip
\noindent
Finally, suppose that $z\in \Lambda(\Gamma)$ is an endpoint of $E_0$.
We claim that $z\in \Lambda (\Phi )$. To prove this, choose
$x\in {\bf B}^n$ which is outside all $\gamma E^0$, $\gamma \in \Gamma$.
Thus there exist $\gamma _i\in \Gamma$ such that $\gamma_i(x)\to z$.
The fact that $x\not \in \bigcup _{\gamma \in \Gamma }\gamma E^0$ implies
$\gamma_i(x)\not \in E^0$. Hence, if $R_i$ is the hyperbolic ray with
endpoints $z$ and $\gamma_i(x)$, then there is
$y_i\in R_i\cap \partial E^0$, since $z$ is an endpoint of $E_0$.
Finally, since $\gamma_i(x) \to z$, we see that $y_i\to z$.
Now fix an $x_0\in \partial E^0$. Since
$\partial E^0/\Gamma =\partial E^0/\Phi $ is compact (because $E^0$ and
$\partial E^0$ are precisely invariant under $\Phi$ in $\Gamma$),
there exist $M>0$ and $\varphi_i\in \Phi $ such that
${\rm d}(\varphi_i(x_0),y_i)\le M$. It follows that $\varphi_i(x_0)\to z$
and hence $z\in \Lambda(\Phi )$. Note that this argument does not 
require that  $E$ is bounded.
\end{proof}

\medskip
\noindent
We now show that the set of endpoints of disjoint ends are disjoint.

\begin{lemma} Let $\Gamma$ be a Kleinian group acting
on ${\bf B}^n$.  If $E^0_1$ and $E^0_2$ are disjoint ends of $\Gamma$, then
their sets of endpoints are disjoint.
Furthermore, if $E^0_2$ is bounded, then the set of endpoints
of $E^0_1$ and the limit set $\Lambda(\Phi_2)$,
$\Phi_2 := \Gamma_{E^0_2}$, are also disjoint.
\label{end-limit-sets-disjoint}
\end{lemma}

\begin{proof}
If $x$ is an endpoint of $E^0_1$, then $x$ is the endpoint of a ray $R$
contained in $E^0_1$ such that ${\rm d}(z,\partial E^0_1) \to \infty $ as
$z\to x$ on $R$. This can be true for at most one end. Therefore, the
sets of endpoints of $E^0_1$ and $E^0_2$ are disjoint, and thus
$\Lambda_e(E^0_1) \cap \Lambda_e(E^0_2) = \emptyset$.
Now, if $E^0_2$ is bounded, then by Lemma~\ref{tukia-lemma-1}
we know that $\Lambda (E^0_2) = \Lambda_c(\Phi_2) \cup \Lambda_e(E^0_2)$.
Clearly, an endpoint of $E^0_1$ cannot be in $\Lambda_c(\Phi_2)$
and the lemma follows.
\end{proof}

\medskip
\noindent
Let ${\cal E}=\{E_1,..,E_n\}$ be a finite collection of ends of
$M_\Gamma$.  For any end $E$ of $M_\Gamma$, let 
$\widetilde{E}=E\cup \partial_\infty E$.  
We say that the collection ${\cal E}$
forms a {\it complete collection of ends} of $M_\Gamma$ if
$E_i$ and $E_j$
are disjoint for $i\ne j$ and if
$\overline{M_\Gamma} \setminus \left( \cup_{E\in {\cal E}}
\widetilde{E}\right)$ is compact.  If ${\cal E}$ forms a complete
collection of ends of $M_\Gamma$ and if in addition each $E_i$ is a
bounded end of $M_\Gamma$, we say that ${\cal E}$ forms a {\em complete
collection of bounded ends} of $M_\Gamma$.

\medskip
\noindent
Let ${\cal E}$ be a complete collection of ends for $M_\Gamma$.  Let
\[
{\cal F}=\{F^0\: |\: F^0\mbox{ is a component of }\pi ^{-1}(E)
\mbox{ for some } E\in {\cal E}\};
\]
we say that ${\cal F}$ is a {\em complete collection of ends} for
$\Gamma$.  If in addition each $F^0\in {\cal F}$ is a bounded end of
$\Gamma$, we say that ${\cal F}$ is a {\em complete collection of bounded
ends} of $\Gamma$. If the [bounded] (pairwise disjoint) ends $F^0_1,...,F^0_k$ of $\Gamma$
generate ${\cal F}$, so that ${\cal F}=\bigcup _i\Gamma F^0_i$, we also
say that $\{F^0_1,...,F^0_k\}$ forms a complete collection of [bounded]
ends for $\Gamma$.

\begin{lemma}
Let $\Gamma$ be a Kleinian group acting on ${\bf B}^n$.
If ${\cal F}$ is a complete collection of ends for $\Gamma$, then
$\Lambda(\Gamma)$ is the disjoint union of $\Lambda _c(\Gamma)$ and of the
end limit point sets $\Lambda _e(F^0)$, $F^0\in {\cal F}$.
\label{tukia-lemma-2}
\end{lemma}

\medskip
\noindent
The proof Lemma \ref{tukia-lemma-2} is similar to that of Lemma
\ref{tukia-lemma-1}.

\section{Conformal measures associated to ends}
\label{measures}

\medskip
\noindent
Let $\Gamma$ be a Kleinian group.  A measure $m$ on a $\Gamma$-invariant 
subset $E$ of $\overline{{\bf B}^n}$ is called an {\em $\alpha$-conformal
measure for $\Gamma$} (or alternatively {\em a conformal measure of dimension $\alpha$
for $\Gamma$}) if all Borel
subsets of  $E$  are measurable and if for any measurable subset  $A$ of
$E$  and for every $\gamma \in \Gamma$, we have
\begin{equation}
\label{conformal}
m (\gamma(A))=\int_{A}|\gamma'|^{\alpha} dm.
\end{equation}
Here, $|\gamma'(\xi)|$ is the operator norm of the derivative of
$\gamma$ at $\xi$, and $\alpha$ is some non-negative number.
Sometimes, when it is clear from the context what group $\Gamma$ and
what number $\alpha$ are meant, we shall just call a measure  
conformal whenever the condition above is satisfied.
Note that we explicitely allow conformal measures to be supported by 
$\overline{{\bf B}^n}$, which stands in contrast to the usual understanding 
that such measures are defined only in the limit set of a group.
Patterson \cite{patterson}, \cite{patterson4} has given the construction
of a probability measure $m$ on $\Lambda(\Gamma)$ which is
$\delta(\Gamma)$-conformal, where $\delta(\Gamma)$ is the critical
exponent of $\Gamma$.
In general, a conformal measure $m$ of dimension $\alpha$ is not unique, 
which is one of the motivations
for this work. However, such an $m$ is unique for groups of divergence type
when $\alpha = \delta(\Gamma)$, see Sullivan \cite{sullivan1},
Roblin \cite{roblin}. There are also other cases when such an $m$ is unique,
see for instance Sullivan \cite{sullivan-growth}.

\medskip
\noindent While Patterson's construction gives 
$\delta(\Gamma)$-conformal measures on $\Lambda(\Gamma)$, it is 
also possible under certain circumstances to construct 
$\alpha$-conformal measures on $\Lambda(\Gamma)$ for $\alpha \geq 
\delta(\Gamma)$.  In fact, an $\alpha$-conformal measure for 
$\Gamma$ can only exist if $\alpha \geq \delta$.  See Sullivan 
(Theorem ($2.19$) of \cite{sullivan-positivity}), and also the 
discussion in Nicholls \cite{nicholls}, Chapter~4.  Thus, the 
critical exponent of $\Gamma$ can be defined as the infimum of 
all numbers $\alpha$ for which there exists an $\alpha$-conformal 
measure on $\Lambda(\Gamma)$.  It is well-known that if the 
Poincar\'e series converges at $\alpha$, then any 
$\alpha$-conformal measure gives zero measure to the conical 
limit set (see for instance \cite{nicholls}, Theorem 4.4.1). A 
construction of such $\alpha$-conformal measures when the 
Poincar\'e series converges at $\alpha$ is given in Theorem 
\ref{theorem-conformal-existence} below. Note that this 
construction is different in nature to Patterson's construction 
\cite{patterson} and in some cases \cite{fatu} it gives different 
measures than Patterson's method.
Both methods work if $\alpha$ is the exponent of convergence and the
Poincar\'e series converges at this exponent; in this case the
non-compactness condition of the theorem is automatically met since
otherwise all limit points are conical limit points and hence the
Poincar\'e series diverges at the exponent of convergence.

\begin{theorem} 
Let $\Gamma$ be a Kleinian group acting on ${\bf B}^n$.  If 
$({\bf B}^n \cup\Omega(\Gamma))/\Gamma$ is not compact, and if the 
Poincar\'e series for $\Gamma$ converges at the exponent $\alpha$, 
then there exists a non-trivial conformal measure of dimension 
$\alpha$ for $\Gamma$ on $\Lambda(\Gamma)$.
\label{theorem-conformal-existence}
\end{theorem}

\begin{proof} 
Since $({\bf B}^n\cup\Omega(\Gamma))/\Gamma$ is not compact, there exists
a sequence of points $\{ z_i\}$ of ${\bf B}^n$ such that the orbits
$\Gamma z_i$ converge to $\Lambda(\Gamma)$ in the Hausdorff metric on
closed subsets of $\overline{{\bf B}^n}$. (In order to see this,
one only needs to choose the $z_i$ so that for any set of the form $\Gamma C$,
where $C\subset {\bf B}^n \cup \Omega(\Gamma)$ is compact, there exists 
$I=I_C$ so that $z_i\not\in \Gamma C$ for $i> I$.  The existence of such
$z_i$ is guaranteed by the assumption that 
$({\bf B}^n\cup\Omega(\Gamma))/\Gamma$ is not compact.) Since the 
Poincar\'e series converges, there exists an atomic $\alpha$-conformal
measure $\mu_i$ for $\Gamma$ on $\Gamma z_i$ of total mass $1$. A
subsequence has a weak limit and this limit is the desired non-trivial
conformal measure on $\Lambda(\Gamma)$.
\end{proof}

\medskip
\noindent
Let now  $E$  be an end of $\Gamma$ in ${\bf B}^n$. Choose the points
$z_i$  in the proof of Theorem~\ref{theorem-conformal-existence}
to be points of  $E$.  We will show that we
obtain a conformal measure supported by the endpoints of $E$.
In order to prove this, we need some estimates on measures of shadows of
hyperbolic balls when viewed from the origin $0\in {\bf B}^n$. 
The next lemma is basically one half of Sullivan's shadow lemma 
(we only need the estimate in one direction), but we sharpen the
statement in the sense that
the constants in the lemma can be chosen not to depend on the measure
$m$ if  $m$ has total mass 1.  We also show that the constant does not
change under conjugation.

\medskip
\noindent
We need only to adapt an argument of Tukia \cite{tukia1} to the present
situation. The half-space model was considered in \cite{tukia1}, as
calculations were easier due to the fact that Euclidean similarities
preserving the half-space are hyperbolic isometries.

\medskip
\noindent
If  $\gamma$ is a M\"obius transformation and  $m$  is a conformal
measure of dimension $\alpha$ on the closed ball $\overline{{\bf B}^n}$,
then we can define a measure  $m_\gamma$, the
{\it image measure} of $\gamma$, by
\begin{equation}
m_\gamma(\gamma (A))=\int_A|\gamma'|^\alpha dm.
\label{image-measure}
\end{equation}
Thus  $m$ is a conformal measure of dimension $\alpha$ for $\Phi$ if
and only if  $m_\gamma=m$ for every $\gamma\in\Phi$.
Considerations involving the Radon-Nikodym derivative show that
$m_{\gamma_1\gamma_2}=(m_{\gamma_2})_{\gamma_1}$ for any two
M\"obius transformations $\gamma_1$ and $\gamma_2$. 
It follows that if $m$ is a conformal measure of dimension $\alpha$ for $\Phi$, then $m_\gamma$ 
is a conformal measure of dimension $\alpha$ for $\gamma \Phi \gamma^{-1}$ 
and is supported by the set $\gamma A$ if  $m$ is supported by $A$.

\medskip
\noindent
We let  $B(z,r)$ be the Euclidean $n$-ball of radius $r$ centred at
$z$.

\begin{lemma}
Let $\Phi$ be a Kleinian group.
Let $C\subset {\bf B}^n$ be compact and fix $k>0$.
Then there exists $M>1$ such that the following holds:  Let $m$  be a
conformal measure of dimension $\alpha$ for $\Phi$  on $\overline{{\bf B}^n}$ of total mass  1.
Let $\gamma$  be a M\"obius transformation
and define $m_\gamma$ as in (\ref{image-measure}).
Consider a point $z\in {\bf S}^{n-1}$ such that for some
$0 \leq t < 1$ we have $(1-t)z\in \gamma (\Phi C)$. Then
\[
m_\gamma (B(z,kt)) \le M t^\alpha.
\]
\label{tukia-lemma-4}
\end{lemma}

\begin{proof}
This is basically Lemma 2C of \cite{tukia1}. In the formulation 
of \cite{tukia1}, we did not claim that the lemma was valid for 
any measure of total mass 1 but rather fixed the conformal 
measure and then found the constants.  However, we need the lemma 
only to have the upper estimate and we need only to set 
$\nu(\overline{{\bf B}^n})=1$ on the second line of p. 247 of the 
proof of Lemma 2C in \cite{tukia1} in order to see that the 
constant in the upper estimate does not depend on the measure if 
the total mass is 1.

\medskip
\noindent
After this observation we transform Lemma 2C of \cite{tukia1} to ${\bf B}^n$
by means of the stereographic projection.  We easily obtain that
Lemma~\ref{tukia-lemma-4} is true if  $z=-e_n=(0,...,0,-1)$
which corresponds to  0 under the stereographic projection.
Other points are obtained by means of an auxiliary
rotation which transforms the point to 0.  Since $|\gamma|'=1$ for a
rotation we obtain our claim in view of the conjugacy invariance of
Lemma 2C of \cite{tukia1}.
\end{proof}

\medskip
\noindent
The following lemma is a direct consequence of Lemma~\ref{tukia-lemma-4}
and is in fact the statement which shall be used in the proof of 
Theorem~\ref{main-theorem}.  Except for conjugacy by a M\"obius
transformation and independence of the measure, assumed to be a
probability measure, it is a direct consequence of Sullivan's shadow lemma.
Here, $S_r(y) \subset \overline{{\bf B}^n}$
denotes the shadow in $\overline{{\bf B}^n}$ from
the origin of the open hyperbolic ball $D(y,r)$ of radius $r>0$
and center $y \in {\bf B}^n$.  Thus a point  $w \in \overline{{\bf B}^n}$
is in $S_r(y)$ if and only if the hyperbolic ray from $0$ to $w$
intersects $D(y,r)$. Note that our definition of the shadow of a 
hyperbolic ball is slightly different from the usual one, where only the 
part contained in ${\bf S}^{n-1}$ is considered.
We apply Lemma~\ref{tukia-lemma-4} in the case that the
compact set $C$ is a one-point set $\{ y \}$.

\begin{lemma}
\label{conjugation-shadowlemma}
Let $\Phi$ be a non-elementary Kleinian group, let $y\in {\bf B}^n$, and
let $r>0$ and $\alpha>0$ be positive constants. Then there is $c>0$ such
that if $m$ is a conformal measure of dimension $\alpha$ for $\Phi$ 
on $\overline{{\bf B}^n}$ of total mass 1, then the following is true.
Let $\gamma$ be an arbitrary M\"obius transformation, and define
the conformal measure $m_\gamma$ for the Kleinian group
$\gamma \Phi \gamma^{-1}$ as in (\ref{image-measure}). Then,
\begin{equation}
\label{sizeofshadows}
m_\gamma(S_r(\gamma (z)))
\; \leq \;
c \; \exp(-\alpha \, {\rm d}(0, \gamma (z)))
\end{equation}
for all $z$ in the orbit  $\Gamma y$.
\end{lemma}

\medskip
\noindent
Armed with these estimates, we can derive some results on the
distribution of mass and extension of conformal measures. We start with
the refinement of Theorem~\ref{theorem-conformal-existence}
which says that we can find a measure
supported by the end limit points.

\begin{theorem}
Let the situation be as in Theorem~\ref{theorem-conformal-existence}  and
let  $E$ be an end of $\Gamma$ in ${\bf B}^n$ such that the end limit 
point set $\Lambda_e(E)\neq\emptyset$.  Then there is an 
$\alpha$-conformal measure $m$ for $\Gamma$  such that $m$
is supported by the end limit points of $E$  and of the ends of
$\Gamma$ equivalent to  $E$ under $\Gamma$. 

\medskip
\noindent
In particular, if the Poincar\'e series for $\Gamma_E$ converges at
exponent $\alpha$ (even if the Poincar\'e series for $\Gamma$ diverges),
there is a non-trivial
conformal measure of dimension $\alpha$ for $\Gamma_E$ supported by
$\Lambda_e(E)$.
\label{conformal-end-measure}
\end{theorem}

\noindent 
{\bf Remark.}
If $\Lambda_e(E)=\emptyset$, then our method still constructs a measure on 
the endpoint set of  $E$  and of the ends equivalent to $E$ under 
$\Gamma$.  However, the set of endpoints which are not end limit points is open and hence the 
$(n-1)$-dimensional Hausdorff measure is a conformal measure on this set.

\medskip
\noindent
\begin{proof}
Pick $y\in\Lambda_e(E)$ and let  $R$  be a hyperbolic ray with endpoint 
$y$.
Pick points $z_i\in R$ such that $z_i\to y$ as $i\to\infty$.  Since  $y$ 
is a limit point of $\Gamma$, it follows that the $z_i$ exit any $\Gamma C$, where $C\subset {\bf B}^n\cup \Omega (\Gamma )$
is compact. Since the Poincar\'e series converges, there is a conformal
measure $m_i$ on $\Gamma z_i$ of total mass 1.
We can assume that the sequence formed by the $m_i$ has a
weak limit $m$.  Since  the 
$\Gamma z_i$ exit any compact subset of ${\bf B}^n\cup\Omega(\Gamma)$, it 
follows that $m$ is supported by $\Lambda(\Gamma)$.  
We show that $m$ is supported by $\Lambda _e(E)$ and the endpoints of the ends equivalent to $E$ under $\Gamma$.

\medskip
\noindent
Let $E_\Gamma =E/\Gamma $ and let $X=\partial E_\Gamma $.
We can assume that $M_\Gamma \setminus X$ has only finitely many 
components. For instance, we can cover $X$ by a finite number of 
hyperbolic balls $B_i$ and replace $E_\Gamma$ by the component of 
$M_\Gamma \setminus (\cup_i B_i)$ containing the original $E_\Gamma$.
Let $F_1,...,F_q$ be the components of $M_\Gamma \setminus X$ whose closures in $M_\Gamma$ are non-compact and which are distinct from  $E_\Gamma$.
Then each $F_i$ is an end of  $M_\Gamma $ and
$\{E_\Gamma,F_1,...,F_q\}$ is a complete collection of ends for $M_\Gamma$.
In view of Lemma \ref{tukia-lemma-2}, it suffices to show that if
$F$ is a lift of some $F_i$ to ${\bf B}^n$, then $m(\Lambda_e(F))=0$.

\medskip
\noindent
We will define for each $0<r<1$ a set $U_r$ so that $U_r$ will be
a neighborhood of $\Lambda _e(F)$ in $\overline{{\bf B}^n}$ and that
$m_i(U_r)\le c_r$ where $c_r\to 0$. This is the basic reason why
$m(\Lambda _e(F))=0$, and the precise argument is given below.

\medskip
\noindent 
Let $H_r=\{z:1-r \le |z|\le 1\}\subset \overline{{\bf 
B}^n}$. We can assume that $0\in E$. 
Let $U_r$ be the union of all geodesic rays 
$R_a=\{ta:1-r < t \leq 1 \}$, $a\in {\bf S}^{n-1}$,
with the property that there exists $1-r < t < 1$
such that $ta \in F$. Clearly, $U_r$ is a neighborhood of 
$\Lambda _e(F)$ in $\overline{{\bf B}^n}$.
We fix a point $z_0\in \partial F$ and a number $R>0$ such that 
$\bigcup _{\gamma \in \Gamma }D(\gamma (z_0),R)\supset \bigcup 
_{\gamma \in \Gamma } \gamma (\partial F)$ where $D(z,R)$ is the 
open hyperbolic ball with center $z$ and radius $R$. Let 
$S_\gamma $ be the shadow of $D(\gamma (z_0),R)$ from $0$, so 
that $S_\gamma $ contains all the points $w\in \overline{{\bf 
B}^n}$ such that the hyperbolic line segment or ray with 
endpoints 0 and $w$ intersects $D(\gamma (z_0),R)$.
Let $V_r$ be the union of all shadows $S_\gamma$, 
$\gamma \in \Gamma$, such that $D(\gamma (z_0),R)$ intersects
$H_r$. It is not difficult to see that 
$\Gamma z_0 \cap U_r \subset V_r$.

\medskip
\noindent
Next, we apply Lemma  \ref{conjugation-shadowlemma} to
the measures $m_i$ whose limit is $m$.
Thus there exists a constant $c$ such that
$m_i(S_\gamma )\le c\, \exp(-\alpha {\rm d}(0,\gamma (z_0)))$
regardless of $i$ and therefore
\[
\sum m_i(S_\gamma )\le c \sum \exp(-\alpha {\rm d}(0,\gamma (z_0))) =: c_r
\]
where both sums are restricted to elements $\gamma\in \Gamma $ such
that $\gamma (z_0)\in H_r$. If $r\to 1$, the right hand side tends to 
zero and thus $c_r \to 0$ as claimed.
Since the mass of $m_i$ in $U_r$ is contained
in the shadows $S_\gamma $, it follows that $c_r$ is
indeed an upper bound for $m_i(U_r)$.

\medskip
\noindent
To see that $m(\Lambda _e(F))=0$, let $\Lambda _p$
be the set of points $z\in {\bf S}^{n-1}$
such that the line segment $tz$, $t\in [1-1/p,1)$, is
contained in $F\cup \partial F$. Note that each $\Lambda_p$ 
is a closed set and $U_r$ is a
neighborhood of $\Lambda _p$ for every $0<r<1$. Since $\Lambda _p$ is
closed, the inequalities $m_i(U_r)\le c_r$ 
imply that $m(\Lambda _p)\le c_r$ for all $r$
and hence $m(\Lambda _p)=0$. Since $\Lambda _e(E)$ is
contained in the union of the $\Lambda _p$, it
follows that  $m(\Lambda _e(E))=0$.

\medskip
\noindent
To see the last paragraph, we only need to observe that $E$ is an end
for $\Gamma_E$ as well, and that the set of end limit points is the same
whether we regard  $E$  as an end of $\Gamma$ or of $\Gamma_E$ 
(see Lemma \ref{tukia-lemma-1}).

\end{proof}

\medskip
\noindent
The above theorem (and the remark following it) 
asserts that we always have a conformal measure for
$\Gamma$ supported by the endpoints of an end  $E$
and the ends equivalent to $E$ under $\Gamma$. 
Conversely, suppose that we have a conformal
measure $m$ for $\Gamma_E$  supported by $\Lambda_e(E)$.  We
might ask whether it is possible to extend $m$ to a conformal measure for
the whole group $\Gamma$. That this is possible is shown in the next proposition.

\begin{proposition}
Let $\Gamma$ be a Kleinian group whose Poincar\'e series converges at 
the exponent $\alpha$, and let $E$ be an end of  $\Gamma$. Let
$m$ be an $\alpha$-conformal measure for $\Gamma_E$
supported by the set of endpoints of $E$.  Then there
is a unique extension of  $m$ to a conformal measure of  $\Gamma$
supported by the endpoints of  $E$  and of the ends equivalent to $E$ 
under $\Gamma$.
\label{extending-conformal-measure}
\end{proposition}

\begin{proof}
Choose representatives $\gamma_i$, $i \in {\bf N}$, from the cosets  
$\Gamma/\Gamma_E$. Let $\gamma_0$ be the identity.  
Thus  $E_i=\gamma_iE$ is an end distinct from  $E=E_0$  if $i\neq 0$.  
Let $\Lambda_i$ be the set of endpoints of $E_i$. 
Thus, the sets $\Lambda_i$, $i \in {\bf N}$, 
form a family of pairwise disjoint sets.
If we can extend $m$ to a conformal measure for $\Gamma$ on 
$\bigcup_i\Lambda_i$, then the restriction of $m$ to $\Lambda_i$ must 
be the image measure $m_{\gamma_i}$ of (\ref{image-measure}). This 
proves the uniqueness of the extension.
Now, the rule $m_{\gamma_1\gamma_2}=(m_{\gamma_2})_{\gamma_1}$, 
and the property that $m_\gamma=m$  if $\gamma\in\Gamma_E$ together imply that
$m_{\gamma_i}$ is independent of the choice of the representative, and
that it satisfies the transformation rule for conformal measures, or
equivalently, that  $m_\gamma=m$ for $\gamma\in\Gamma$.
These properties also imply that setting  $m:=m_{\gamma_i}$ on $\Lambda_i$,
we obtain a measure on
$\bigcup_i \Lambda_i$ which satisfies the transformation rule for conformal
measures.  Thus it is a conformal measure for $\Gamma$ if it is finite.

\medskip
\noindent
We will now prove the finiteness of $m$ by relating $m(\bigcup \Lambda_i)$
to the sum of the Poincar\'e series at $\alpha$, which is assumed to be finite.  
It suffices to prove that $\sum_{i\neq 0}m(\Lambda_i)<\infty$.

\medskip
\noindent
By conjugation with suitable M\"obius transformations
we can assume that $0\in E$. Next, we pick $y\in \partial E$.  
Since  $\partial E/\Gamma$ is compact, there is a number  $r>0$
such that the shadows from $0$ of the hyperbolic balls  
$D(\gamma(y),r)$, $\gamma\in\Gamma$,
cover $\bigcup_i \partial E_i$.  Suppose that  $z\in \Lambda_i$ for some $i\neq 0$.
Then the hyperbolic ray with endpoints  0  and $z$ intersects
$\partial E_i$ and hence this ray intersects also some  $D(\gamma(y),r)$
so that $z\in S_r(g\gamma_i(y))$ for some $g\in
\Gamma_{E_i}=\gamma\Gamma_E\gamma^{-1}$.
Lemma~\ref{conjugation-shadowlemma} then implies the
existence of a constant $C>0$  such that
\[
m \bigg(\bigcup_{i\neq 0} \Lambda_i\bigg)
\leq m \bigg( \bigcup_{\gamma\in\Gamma}\,S_r(\gamma(y)) \bigg)
\leq \sum_{\gamma\in\Gamma}C \exp(-\alpha {\rm d}(0,\gamma(y))) <\infty.
\]
\end{proof}

\medskip
\noindent
{\bf Remark.}
We have formulated Theorem~\ref{conformal-end-measure} and
Proposition~\ref{extending-conformal-measure} for the case at hand 
so that we take the weak limit of atomic measures supported by
an orbit or extend measures supported by the end limit point set.
A more general formulation would be as follows.
Let $E$ be an end of $\Gamma$ in ${\bf B}^n$ and set
$\widetilde{E} = 
(\overline{E}\cap({{\bf B}^n}\cup\Omega(\Gamma)))\cup\Lambda_e(E)$.
The formulation of 
Proposition~\ref{extending-conformal-measure} 
would be that an $\alpha$-conformal measure $\mu$ for $\Gamma_E$
which is supported by $\widetilde{E}$ can be extended to an 
$\alpha$-conformal measure for $\Gamma$
supported by $\bigcup_{\gamma\in\Gamma} \widetilde{E}$.
Theorem~\ref{conformal-end-measure} is formulated in the clearest way for
conformal measures for $\Gamma_E$ (thus $\Gamma=\Gamma_E$). In this case,
if the measures $\mu_i$ are supported by $\widetilde{E}$ and have total
mass 1, then also their weak limit $\mu$ is supported by $\widetilde{E}$.

\medskip
\noindent
Suppose that ${\cal E}=\{E_1,...,E_q\}$ is a complete collection of ends for
$\Gamma $ where each  $E_i$ is an end of $\Gamma$ in ${\bf B}^n$ and where the $E_i$ are pairwise disjoint.  Thus, setting $E_{i\Gamma }=E_i/\Gamma $, it follows that
$M_\Gamma \setminus (\bigcup _{i\leq q}E_{i\Gamma })$ is compact.
Let $\Gamma_i := \Gamma_{E_i}$. Fix $\alpha >0$
such that the Poincar\'e series for $\Gamma $ converges
with exponent $\alpha $. We suppose that  $E_i$, $i\leq p$, are the ends  
such that  $\Lambda_e(E_i)\neq\emptyset$ and denote
\medskip

${\cal M}=$ the family of $\alpha$-conformal measures 
for $\Gamma $ on $\Lambda (\Gamma )$.

${\cal M}_i=$ the family of $\alpha$-conformal measures 
for $\Gamma _i$ on $\Lambda _e(E_i)$, $i\leq p$.

\medskip
\noindent
We can now prove our main theorem.  The natural situation for us is that 
the measures live on $\Lambda(\Gamma)$ but everything remains valid if ${\cal 
M}$ is the set of all conformal measures on ${\bf S}^{n-1}$ and  
${\cal M}_i$  is the set of conformal measures on endpoints of $E_i$
as $i$ varies from 1 to $q$.

\begin{theorem}
Let $\mu \in {\cal M}$ be a conformal measure for $\Gamma $  on $\Lambda(\Gamma)$, and let
$\mu_i\in {\cal M}_i$ be the restriction of  $\mu$ to $\Lambda_e(E_i)$.
Then
\begin{equation}
   \mu =\mu_1^*+...+\mu_p^*
\label{decomposition-of-measure}
\end{equation}
where $\mu _i^*$ is the unique extension of $\mu _i$ to a conformal 
measure on $\bigcup _{\gamma \in \Gamma }\Lambda _e(E_i)$ given by 
Proposition~\ref{extending-conformal-measure}.
The measures $\mu_i^*$ and $\mu_j^*$, $i\neq j$, are mutually singular
if non-trivial.

\medskip
\noindent
If $\mu_i\in {\cal M}_i$, then
(\ref{decomposition-of-measure}) defines a conformal
measure $\mu $ for $\Gamma $ on $\Lambda (\Gamma )$.
For each  $i\leq p$ there is a non-trivial measure $\mu_i\in{\cal M}$ and thus,
if there are $p$ ends in ${\cal E}$  such that 
$\Lambda_e(E_i)\neq\emptyset$, 
then there are at least $p$ mutually singular non-trivial conformal
measures of dimension $\alpha $  for $\Gamma $ on $\Lambda (\Gamma )$.
\label{main-theorem}
\end{theorem}

\begin{proof}
By the uniqueness of the extension of
Proposition~\ref{extending-conformal-measure},
$\mu _i^*$ and $\mu $ coincide on
$E_i^*=\bigcup _{\gamma \in \Gamma }\gamma \Lambda _e(E_i)$.
Every limit point of $\Gamma $ is either in some
$E_i^*$ or is a conical limit point of $\Gamma $
(Lemma~\ref{tukia-lemma-2}).
Since the Poincar\'e series converges, conical limit points have
zero measure and hence (\ref{decomposition-of-measure}) is true.
Since $\mu_i^*$ is supported by $E_i^\ast$,  and the $E_i^\ast$ are disjoint,
it follows that $\mu_i$ and $\mu_j$ are mutually singular if non-trivial.

\medskip
\noindent
To obtain the last paragraph, we note that 
Theorem~\ref{conformal-end-measure} gives the non-trivial measure 
$\mu_i$ on $\Lambda_e(E_i)$ which can be extended to the conformal 
measure $\mu_i^\ast$ by Proposition~\ref{extending-conformal-measure}.
Other points of the last paragraph are obvious.
\end{proof}

\medskip
\noindent
{\bf Bounded ends.}
Theorem~\ref{main-theorem} says that it is possible to write a
conformal measure $\mu $ for $\Gamma $ on $\Lambda (\Gamma )$
as a sum of extended measures, which are obtained from 
measures supported by the end limit point sets of ends. This is a way to
decompose $\mu $ to simpler measures. The decomposition of
(\ref{decomposition-of-measure}) is not ideal,
since $\mu_i$ is supported by the end
limit set of $\Gamma_i=\Gamma_{E_i}$. It would be
better if we could allow $\mu_i$ to be 
supported by $\Lambda (\Gamma_i)$ since, for
instance, the Patterson-Sullivan construction of a conformal measure
gives measures on $\Lambda (\Gamma_i)$ and it seems that these measures
are not necessarily supported by the end limit points.

\medskip
\noindent
This problem disappears if the end is bounded since then every
$z\in \Lambda (\Gamma _i)$ is either an end limit point
or a conical limit point by Lemma \ref{tukia-lemma-1} and the
measure of the conical limit set vanishes in the case of convergence of the
Poincar\'e series. Thus, in this case we can
replace $\Lambda_e(E_i)$ by $\Lambda (\Gamma_i)$.

\medskip
\noindent
We combine our main theorems applied to bounded ends to

\begin{theorem}
Let $E$ be a bounded end of $\Gamma$ and suppose that the Poincar\'e
series for $\Gamma$ converges at the exponent $\alpha$.
If $\mu $ is a $\alpha$-conformal measure for $\Gamma _E$ on
$\Lambda (\Gamma _E)$, then $\mu$ has a unique extension, 
denoted by $\mu^*$, to a conformal measure for $\Gamma $ such
that the extension is supported by the subset
$\,\bigcup_{\gamma \in \Gamma }\gamma \Lambda (\Gamma _E)$ of
$\Lambda (\Gamma )$.

\medskip
\noindent
Let $\{E_1,...,E_p\}$ be a complete collection of ends for $\Gamma $ 
such that each $E_i$ is a bounded end in ${\bf B}^n$.
If $\mu $ is a conformal measure of dimension $\alpha $ on
$\Lambda (\Gamma )$, then $\mu =\mu _1^*+...+\mu _p^*$
where $\mu _i$ is the restriction of $\mu$ to $\Lambda (\Gamma _{E_i})$,
and the $\mu _j^*$, $i\neq j$, are mutually singular if they are non-trivial.
\label{bounded-end-theorem}
\end{theorem}

\medskip
\noindent
{\bf Applications.}
We have the following consequences of the measure extension
Proposition~\ref{extending-conformal-measure}.
The first theorem is proved using the fact, discussed in Example 2 in Section \ref{ends-and-endgroups}, that if $v$ is a parabolic fixed
point of full rank, then its stabilizer $\Gamma_v$ is a bounded end group.
Thus the atomic measure which gives mass 1 to $v$ is conformal for
$\Gamma_v$ and hence can be extended to a conformal measure for $\Gamma$.

\begin{theorem}
Let the Kleinian group $\Gamma$ be as in
Theorem~\ref{main-theorem}, and let $v\in {\bf S}^{n-1}$ be a parabolic 
fixed point such that $\Gamma_v$ has rank $n-1$. 
Then there is a non-trivial atomic conformal measure for $\Gamma$ 
supported by the orbit $\Gamma v$ of $v$. If $g_i$, $i\in I$, are
representatives of the cosets of $\Gamma/\Gamma_v$, then
\[
\sum _{i\in I}|g'_i(v)|^\alpha < \infty .
\]
\label{consequence-1}
\end{theorem}

\begin{theorem}
Let  $\Gamma$ be a Kleinian group and let $E$ be an end of $\Gamma$ in
${\bf B}^n$.
Let $m$ be a conformal measure
for $\Gamma_E$ on $\Lambda(\Gamma_E)$. If $g_i$, $i\in I$,
are representatives of the cosets of $\Gamma/\Gamma_E$, then
\begin{equation}
\sum _{i\in I}|g'_i(x)|^\alpha < \infty
\label{finite}
\end{equation}
for $m$-almost every $x\in \Lambda_e(E)$.  If the end $E$ is bounded,
then (\ref{finite}) is true for $m$-almost every $x\in
\Lambda(\Gamma_E)$.
\label{consequence-2}
\end{theorem}

\medskip
\noindent
{\bf Ends with boundary.}
We have defined ends as subsets $E$ of $M_\Gamma $ (or the lift to
${\bf B}^n$ of such a set) such that $E$ is a non-compact component of
$M_\Gamma \setminus \partial E$ where
the boundary $\partial E$ of $E$ in $M_\Gamma $ is compact.
We could have done this
in $\overline {M}_\Gamma =({\bf B}^n\cup \Omega (\Gamma ))/\Gamma $
so that $E$ would now be a component of
$\overline {M}_\Gamma \setminus \bar \partial E$
where now the boundary $\bar \partial $
is taken in $\overline {M}_\Gamma $. We call such
a set $E$ an {\it end of $M_\Gamma$ with boundary}
(of course, such an end does not necessarily intersect the boundary
of $M_\Gamma )$. Like before, we also call a lift
$\widetilde{E}$ of $E$ to ${\bf B}^n\cup \Omega (\Gamma )$
an end of $\Gamma$ with boundary. Such an end is bounded if
$({\bf B}^n\cup \Omega (\Gamma_{\widetilde{E}} ))\setminus \widetilde{E}$
has compact $\Gamma_E$-quotient. The definition of the endpoints and end
limit points of an end with boundary can be given as before.

\medskip
\noindent
All the earlier results are valid with this definition.
Lemma~\ref{tukia-lemma-1} and its tripartite division of points 
of ${\bf S}^{n-1}$ into ordinary, conical limit
points or end limit points is valid. Similarly, if we define a
complete collection of ends with boundary to be a collection
${\cal E}=\{E_1,...,E_p\}$ such that 
$\overline {M}_\Gamma \setminus (E_1\cup ...\cup E_p)$ is compact, 
or the lift of such a collection, then Lemma~\ref{tukia-lemma-2} is valid,
and every $z\in \Lambda (\Gamma )$ is either a conical limit
point or an endpoint of a lift to ${\bf B}^n$ of some  $E_i$. If an 
endpoint  $z$ of $E_i$ is in $\Lambda(\Gamma)$, then 
$z\in\Lambda(\Gamma_{E_i})$.

\medskip
\noindent
Finally, as in Theorem~\ref{conformal-end-measure} we can find a non-trivial
conformal measure supported by the end limit point set of an end
$E$ in ${\bf B}^n$ and Proposition~\ref{extending-conformal-measure} 
is still valid, allowing to extend a conformal measure for $\Gamma_E$
supported by the end limit points to a conformal measure for $\Gamma$.
Note that in the proof of these statements we now 
should choose  that
$z_0$ and the point 0, from which we take shadows,  are in the
hyperbolic convex hull $H_\Gamma$ of
$\Lambda(\Gamma)$. (The convex hull is discussed in Section \ref{ends-and-endgroups-3dim}.) Thus the hyperbolic line or segment joining 0
and a point of some  $E_i$  or a point of $\Gamma z_0$  lies in 
$H_\Gamma$.  Therefore, in the proofs of
Proposition~\ref{extending-conformal-measure} and 
Theorem~\ref{conformal-end-measure}
we can replace $F$ by $F\cap H_\Gamma$ and $E_i$ by
$E_i\cap H_\Gamma$, both of which can be covered by the 
shadows of the balls $D(\gamma(z_0),R)$ and
$D(\gamma(y),R)$,  $\gamma\in\Gamma$.

\medskip
\noindent
In analogy to Theorem~\ref{main-theorem} we could state that if
${\cal E}=\{E_1,...,E_p\}$ is a complete collection of ends with boundary, then
{\it a conformal measure $\mu $ of $\Gamma $ on $\Lambda(\Gamma)$
of dimension $\alpha $ such that the Poincar\'e series converges
at $\alpha $, admits a unique decomposition
as a sum $\mu =\mu _1^*+....+\mu _p^*$ where $\mu _i$
is a conformal measure supported
by the end limit points of $E_i$ and $\mu _i^*$ is the extension of 
$\mu_i$ to a conformal measure for $\Gamma$ which is given by the 
analogue of Proposition~\ref{extending-conformal-measure} for ends 
with boundary.}
Also, other parts of Theorem~\ref{main-theorem}
are valid for ends with boundary and so is the analogue of
Theorem~\ref{bounded-end-theorem} for bounded ends with boundary.

\section{Ends and end groups of Kleinian groups acting on ${\bf B}^3$}
\label{ends-and-endgroups-3dim}

\medskip
\noindent
We now specialize to the case of a Kleinian group acting on ${\bf B}^3$, and discuss the behavior of ends of hyperbolic $3$-manifolds in the context of the terminology of this paper.  Throughout this section, let $\Gamma$ be a non-elementary, finitely generated, purely loxodromic, torsion-free Kleinian group acting on ${\bf B}^3$.  We further assume that the Kleinian groups considered in this section contain only orientation-preserving isometries of ${\bf B}^3$. The assumptions that $\Gamma$ is purely loxodromic and torsion-free are not essential, but are made for ease of exposition.   

\medskip
\noindent
A Kleinian group $\Gamma$ acting on ${\bf B}^3$ is 
{\em topologically tame} if $M_\Gamma$ is
homeomorphic to the interior of a compact $3$-manifold with (possibly
empty) boundary.  In particular, topologically tame Kleinian groups are
finitely generated.  It is conjectured that all finitely
generated Kleinian groups are topologically tame.  
Agol \cite{agol} and Calegari and Gabai \cite{calegari-gabai} 
have recently and independently announced proofs of this conjecture.   We note that in the case that $\Gamma$ is topologically tame and geometrically infinite, 
it is known that $\delta(\Gamma) = 2$, see Canary \cite{canary1}, 
\cite{canary2}.  Again for ease of exposition, we restrict our attention in this section to topologically tame Kleinian groups.

\medskip
\noindent
Let $E^0$ be an end of $\Gamma$ with associated end group $\Gamma_{E^0}$. Say that $E^0$ is {\it finite} if
$\overline{E^0}\cap {\bf S}^2=\overline{D}$, where $D$ is a component of
$\Omega(\Gamma)$.  Say that $E^0$ is {\it infinite} if 
$\overline{E^0}\cap {\bf S}^2\subset \Lambda(\Gamma)$.   As we will see below, an end of $\Gamma$ is not necessarily either finite or infinite.
By Lemma \ref{tukia-lemma-1}, the end $E^0$ of $\Gamma$ being infinite is
equivalent to the fact that 
$\overline{E^0}\cap {\bf S}^2\subset \Lambda(\Gamma _{E^0})$.  
We note that the ends of $M_\Gamma$ that are
finite by this definition correspond to the ends of $M_\Gamma$ that are
geometrically finite ends by the more standard definition.

\medskip
\noindent
The {\em convex core} $C_\Gamma$ of the hyperbolic $3$-manifold
$M_\Gamma = {\bf B}^3/\Gamma$ is the smallest closed, convex subset of
$M_\Gamma$ whose inclusion into $M_\Gamma$ is a homotopy equivalence. 
Equivalently, $C_\Gamma$ is the quotient under $\Gamma$ of the
{\em convex hull} $H_\Gamma$ of $\Lambda(\Gamma)$ in ${\bf B}^3$, where the convex
hull $H_\Gamma$ of $\Lambda(\Gamma)$ is the smallest closed, convex subset of
${\bf B}^3$ containing all the lines in ${\bf B}^3$ both of whose
endpoints at infinity lie in $\Lambda(\Gamma)$.
A Kleinian group $\Gamma$ is {\em geometrically finite} if a unit
neighborhood of $C_\Gamma$ has finite volume,
and is {\em geometrically infinite} otherwise.   A Kleinian group acting on
${\bf B}^2$ is geometrically finite if and only if it is finitely
generated.  A geometrically finite Kleinian group acting on ${\bf B}^3$ is necessarily finitely generated, but the converse fails.

\medskip
\noindent
A {\em compact core} for $M_\Gamma$ is a compact submanifold $Y$ of $M_\Gamma$ for which the inclusion map $Y \hookrightarrow M_\Gamma$ induces a homotopy equivalence.  It is a theorem of Scott \cite{scott-core} that for any finitely generated Kleinian group $\Gamma$ acting on ${\bf B}^3$ (in fact, for any irreducible, orientable $3$-manifold with finitely generated fundamental
group), there exists a compact core for $M_\Gamma$.  Note that if $\Gamma$
is topologically tame, then $M_\Gamma$ is homeomorphic to the interior of
a compact $3$-manifold $Z$ with boundary, and a compact core for
$M_\Gamma$ can be obtained by removing a collar neighborhood of 
$\partial Z$ from $Z$.  In this case, the connected components of $M_\Gamma \setminus Y$ form a complete collection of ends and so, in the case that $\Gamma$ is topologically tame, there exists a 
complete collection of ends ${\cal E}$ whose elements are in one-to-one
correspondence with the boundary components of $Z$.   

\medskip
\noindent
Let $E$ be an end of $M_\Gamma$.  We say that $E$ is a {\em peripheral end} of $M_\Gamma$ if there exists a compact core $Y$ of $M_\Gamma$ so that $E$ is contained in a component of $M_\Gamma \setminus Y$. Note that there is a great deal of flexibility in the definition of a peripheral end, as we have a great deal of flexibility in choosing the compact core.  (We acknowledge that this is not the most general definition of peripheral end that can be given.  However, given the generality with which we are treating ends, it seems an appropriate definition.)  With this definition, a peripheral end is a neighborhood of an end in the sense of Bonahon \cite{bonahon}.  A peripheral end is finite if it faces a component of $\Omega(\Gamma)/\Gamma$.   If $E$ is a peripheral end of $M_\Gamma$, then we call any component $E^0$ of the lift of $E$ to ${\bf B}^3$ a peripheral end of $\Gamma$. 

\medskip
\noindent
There is one basic class of Kleinian groups that are of interest to us here.   Let $\Gamma$ be a purely loxodromic Kleinian group isomorphic to the fundamental group of a closed orientable surface $S$  of negative Euler characteristic.  By Bonahon's criterion for tameness \cite{bonahon}, each such $M_\Gamma$ is homeomorphic to the interior of a compact $3$-manifold with boundary, and basic $3$-manifold topology, see e.g. Hempel \cite{hempel}, implies that $M_\Gamma = S\times (0,1)$.  

\medskip
\noindent
For any $t$ in $(0,1)$, the surface $S_t := S\times \{ t\}$ is compact and separating, and so the two components $E_{(0,t)} := S\times (0,t)$ and $E_{(t,1)} := S\times (t,1)$ of $M_\Gamma \setminus S_t$ are peripheral ends of $M_\Gamma$.  In fact, if $E$ is any peripheral end of $M_\Gamma$, then, for any $t\in (0,1)$, the symmetric difference of $E$ with one of $E_{(0,t)}$ or $E_{(t,1)}$ has compact closure, and so to describe the peripheral ends of $M_\Gamma$, it suffices to consider $E_{(0,t)}$ and $E_{(t,1)}$.  (To see this, note that if $E$ is any peripheral end of $M_\Gamma$, then $\partial E\cup \partial E_{(0,t)}\cup \partial E_{(t,1)}$ is compact, and so is contained in $S\times [a,b]$ for some closed interval $[a,b]$.  Then, we see that $E \setminus S\times [a,b]$ must coincide with either $E_{(0,t)}\setminus S\times [a,b]$ or $E_{(t,1)} \setminus S\times [a,b]$, say $E_{(0,t)}\setminus S\times [a,b]$.  In particular, the symmetric difference of $E$ and $E_{(0,t)}$ is then contained in $S\times [a,b]$, which is compact.)  In fact, since $M_\Gamma$ is the product $S\times (0,1)$, we see that for any peripheral end $E$ of $M_\Gamma$, the other component of $M_\Gamma \setminus \partial E$  is itself a peripheral end of $M_\Gamma$, which we refer to as the {\em complementary end to $E$ in $M_\Gamma$}.

\medskip
\noindent
Since $S_t$ is incompressible, it lifts to a properly embedded open topological disc $S^0_t$ in ${\bf B}^3$, which is necessarily invariant under $\Gamma$.  The two components $E_{(0,t)}^0$ and $E_{(t,1)}^0$ of ${\bf B}^3 \setminus S_t^0$ are connected, invariant under $\Gamma$, and not equivalent under $\Gamma$.  Again up to taking the symmetric difference with the $\Gamma$-translates of a compact set, these are the only two peripheral ends of $\Gamma$.    Note that since both $E_{(0,t)}^0$ and $E_{(t,1)}^0$ are invariant under $\Gamma$, we have that $\Lambda(\Gamma)\subset \overline{E_{(0,t)}^0}$ and $\Lambda(\Gamma)\subset \overline{E_{(t,1)}^0}$.  

\medskip
\noindent
The domain of discontinuity $\Omega(\Gamma)$ of $\Gamma$ contains at most two components, and all components of $\Omega(\Gamma)$ are invariant under $\Gamma$.   If the peripheral end $E^0$ of $\Gamma$ does not face a component of $\Omega(\Gamma)$, then (up to taking its symmetric difference with the $\Gamma$-translate of a compact set) $E^0$ is one of $E_{(0,t)}^0$ or $E_{(t,1)}^0$ and is contained in the convex hull $H_\Gamma$ of $\Lambda(\Gamma)$ in ${\bf B}^3$; in particular, $\overline{E^0}$ is disjoint from $\Omega(\Gamma)$, as any geodesic ray in ${\bf B}^3$ ending at a point of $\Omega(\Gamma)$ must exit the convex hull $H_\Gamma$ of $\Lambda(\Gamma)$ in ${\bf B}^3$.  Hence, if $E^0$ does not face a component of $\Omega(\Gamma)$, then $E^0$ is an infinite peripheral end.

\medskip
\noindent
To summarize:

\begin{proposition} Let $\Gamma$ be a purely loxodromic Kleinian group isomorphic to the fundamental group of a closed orientable surface $S$ of negative Euler characteristic.  Then, every peripheral end of $\Gamma$ is either finite or infinite.  

\medskip
\noindent
A peripheral end of $\Gamma$ is finite if and only if it faces a component of $\Omega(\Gamma)$, and a peripheral end of $\Gamma$ is infinite if and only if it is contained in the convex hull $H_\Gamma$ of $\Lambda(\Gamma)$ (up to having symmetric difference contained in the $\Gamma$-translates of a compact set).   

\medskip
\noindent
If the peripheral end $E$ of $M_\Gamma$ is finite, then every endpoint of $E$ is a point of $\Omega(\Gamma)$, and if $E$ is infinite, then every endpoint of $E$ is an end limit point of $E$.
\label{finite-infinite-prop}
\end{proposition}

\medskip
\noindent
Consider now the peripheral end $E_{(0,t)}^0$ of $\Gamma$.  Suppose that its complementary end $E_{(t,1)}^0$ faces a component $\Delta_1$ of $\Omega(\Gamma)$.  Set 
\[ F_{(0,t)}^0 = (\overline{{\bf B}^3} \setminus \overline{E_{(0,t)}^0}) \cup\partial E_{(0,1)}^0. \]  
Since 
\[ \overline{{\bf B}^3} \setminus \overline{E_{(0,t)}^0} = E_{(t,1)}^0\cup \Delta_1 \]
and $\partial E_{(0,t)}^0 = \partial E_{(t,1)}^0$, we have that
\[ F_{(0,t)}^0 = E_{(t,1)}^0\cup \partial E_{(t,1)}^0\cup \Delta_1, \]
and so $F_{(0,t)}^0/\Gamma$ is compact.  That is, if $E$ is a peripheral end of $M_\Gamma$ and if its complementary end faces a component of $\Omega(\Gamma)$, then $E$ is bounded.  If its complementary end does not face a component of $\Omega(\Gamma)$, then $E$ is not bounded, as the boundary at infinity of its complementary end is empty.  

\medskip
\noindent
To summarize:

\begin{proposition} Let $\Gamma$ be a purely loxodromic Kleinian group isomorphic to the fundamental group of a closed orientable surface $S$  of negative Euler characteristic.  Then, a peripheral end of $\Gamma$ is bounded if and only if its complementary end is finite.  
\label{bounded-notbounded-prop}
\end{proposition}

\medskip
\noindent
Say that $\Gamma$ is {\em quasifuchsian} if $\Lambda(\Gamma)$ is a Jordan curve and no element of $\Gamma$ interchanges the two components of ${\bf S}^2 \setminus \Lambda(\Gamma)$.  In particular, $\Omega(\Gamma)$ has two components, each invariant under $\Gamma$.  In this case, both peripheral ends of $\Gamma$ are finite and both peripheral ends of $\Gamma$ are bounded.

\medskip
\noindent
Say that $\Gamma$ is {\it degenerate} if $\Omega(\Gamma)$ is non-empty and
consists of one simply connected component.  In this case, $\Gamma$ has one finite peripheral end and one infinite peripheral end.  The peripheral end not facing the one component of $\Omega(\Gamma)$ is infinite, while the peripheral end facing the component of $\Omega(\Gamma)$ is finite.  Moreover,  the infinite peripheral end is bounded, while the finite peripheral end is not bounded.   A degenerate group $\Gamma$ for which there is a global positive
lower bound on the injectivity radius over all of $M_\Gamma$ is sometimes
known as a {\em hyperbolic half-cylinder}; these are the hyperbolic
$3$-manifolds considered by Sullivan in \cite{sullivan-growth}.

\medskip
\noindent
Say that $\Gamma$ is {\em doubly degenerate} if $\Gamma$ if $\Lambda(\Gamma) ={\bf S}^2$.  
Examples of doubly degenerate groups include the fiber covering of a
closed hyperbolic $3$-manifold fibering over ${\bf S}^1$, though there are examples of doubly degenerate groups that are not associated to fibered $3$-manifolds.   In this case, both peripheral ends of $\Gamma$ are infinite, and neither peripheral end of $\Gamma$ is bounded.

\medskip
\noindent
In general, we have the following basic construction of the complete collection of ends of any
non-elementary, finitely generated, purely loxodromic Kleinian group
acting on ${\bf B}^3$.  (We note that the following proposition also holds
for elementary Kleinian groups, but we leave the proof to the interested
reader.)

\begin{proposition} Let $\Gamma$ be a purely loxodromic, topologically tame 
Kleinian group acting on ${\bf B}^3$.  Then $M_\Gamma$ has a complete collection of ends such that
every end $E\in {\cal E}$ is peripheral and every end $E\in {\cal E}$ is either finite or infinite.
\label{tukia-lemma-3}
\end{proposition}

\begin{proof} 
Let $Y$ be a  compact core for $M_\Gamma$.  Note that $\partial Y$ is empty if and only if 
$M_\Gamma$ is itself compact.  In this case, the collection of ends of $M_\Gamma$ 
is empty.

\medskip
\noindent
Suppose now that $M_\Gamma$ has infinite volume.  The ends will be the 
components of the complement of $Y$.   Enumerate the components of 
$\partial Y$ as $\partial Y =\cup_{j=1}^p S_j$, where  $S_j$ faces a component  of $\Omega(\Gamma)/\Gamma$ for $1\le j\le \ell$, where $\ell \le p$.  Let $E_j$ be the
component of $M_\Gamma \setminus Y$ facing $S_j$.  Note that
$S_j =\partial E_j$ is a connected component of $\partial Y$, and hence is a
separating surface in $M_\Gamma$, and that $E_j$ is a non-compact
component of $M_\Gamma \setminus \partial E_j$, and so
${\cal E} =\{ E_1,\ldots, E_p\}$ forms a complete collection of ends for
$M_\Gamma$.  Moreover, by construction, each end in ${\cal E}$ is peripheral.

\medskip
\noindent
It remains only to determine the types of the ends in ${\cal E}$.  The
natural retraction of $\overline{M_\Gamma}$ onto $C_\Gamma$ yields that the
ends $E_1,\ldots, E_\ell$ facing the components of $\Omega(\Gamma)/\Gamma$ are
finite.  The ends $E_{\ell+1},\ldots, E_p$ are essentially contained in the convex
core $C_\Gamma$ of $M_\Gamma$, and hence are infinite.  (By {\em essentially} here, we mean that for each $\ell+1\le j\le p$, there is a peripheral end $E'_j$ of $M_\Gamma$ contained in the convex core $C_\Gamma$ of $M_\Gamma$ for which the symmetric difference of $E_j$ and $E'_j$ has compact closure.)
\end{proof}

\medskip
\noindent
We note here that with the definition given here, there are many ends that are not peripheral, as they are not components of the complement of a compact core.  The following example contains an example of such a non-peripheral end.   Consider a non-elementary, finitely generated, purely loxodromic
Kleinian group $\Gamma$ whose quotient $M_\Gamma$ has the following
structure: Let $N$ be a compact hyperbolizable acylindrical $3$-manifold
with $3$ boundary components $S_1$, $S_2$, and $S_3$.  Put a hyperbolic
structure on the interior ${\rm int}(N) = M_\Gamma$ in such a way that the
inclusion of $\pi_1(S_1)$ into $\Gamma$ is Fuchsian, the inclusion of
$\pi_1(S_2)$ into $\Gamma$ is quasifuchsian, and the inclusion of
$\pi_1(S_3)$ into $\Gamma$ is a degenerate group.  (Such a manifold $N$
can be constructed by taking the complement of a sufficiently complicated
$3$ component graph in ${\bf S}^3$, and such a hyperbolic structure can be
constructed by taking an appropriate limit of geometrically finite
hyperbolic structures on ${\rm int}(N)$.)  Doubling across $S_1$, we
obtain a hyperbolic $3$-manifold $P$ so that $P$ contains a separating
totally geodesic surface $S$ (corresponding to $S_1$).  The two components
of $P \setminus S$, namely the two ends $E_0$ and $E_1$ associated to $S$, 
are isomorphic (by reflection across $S$), and $S =\partial E_k$ is incompressible in $P$.  Neither end is finite nor infinite, since it is not contained in the convex core of $P$, and neither end faces a component of the domain of discontinuity of the Kleinian group uniformizing $P$.  (In fact, each end has subends that do both of these.) In this case, though, both ends are bounded: since the subgroup $\Phi$ of the Kleinian group uniformizing $P$ corresponding to $S$ is Fuchsian, the complementary end of each $E_k$ is half of the Fuchsian manifold $\overline{M_\Phi}$.  (In fact, all that we need is that $M$ be a compact hyperbolizable $3$-manifold containing an embedded incompressible separating surface $S$ that is not homotopic into $\partial M$.  We have just given an explicit construction of such $M$ and $S$.)

\medskip
\noindent
Now, let $\Gamma$ be a purely loxodromic, topologically tame Kleinian group, let $E$ be an end of $M_\Gamma$, and suppose that $\partial E$ is a separating surface in $M_\Gamma$.  (We restrict to the case that $\partial E$ is a surface for ease of exposition.)  Say that $E$ is an {\em incompressible end} if $\partial E$ is incompressible and if the inclusion of $\partial E$ into $E\cup\partial E$ induces an isomorphism of fundamental groups.  We can characterize the incompressible ends of $M_\Gamma$.

\begin{proposition}  Let $\Gamma$ be a purely loxodromic, topologically tame Kleinian group which is not isomorphic to the fundamental group of a closed, orientable surface of negative Euler characteristic, and let $E$ be an incompressible end of $M_\Gamma$ with end group $\Phi$.  Then, either $\Phi$ is quasifuchsian and $E$ is finite, or $\Phi$ is degenerate and $E$ is infinite; in either case, $E$ is bounded.
\label{incompressible-bounded}
\end{proposition}

\begin{proof} Since $\partial E$ is an incompressible surface, we have that $\Phi$ is a purely loxodromic Kleinian group that is isomorphic to the fundamental group of a closed, orientable surface of negative Euler characteristic, and hence is either  quasifuchsian, degenerate, or simply degenerate.  (This fact is standard.  See e.g. Anderson \cite{anderson}.)  Moreover, by the assumption made on $\Gamma$, it must be that $\Phi$ has infinite index in $\Gamma$.  

\medskip
\noindent
This immediately implies that $\Phi$ cannot be doubly degenerate: if $\Phi$ is doubly degenerate, then by the end covering theorem of Canary \cite{canary-covering}, the covering $M_\Phi \rightarrow M_\Gamma$ is one-to-one, a contradiction.  So, $\Phi$ is either quasifuchsian or degenerate.  Note that $E$ can then also be considered an end of $M_\Phi$.  In particular, $E$ is either finite or infinite, since this dichotomy holds for ends of quasifuchsian and degenerate groups, see Proposition \ref{finite-infinite-prop}.

\medskip
\noindent
It remains only to show that $E$ cannot be the finite end of a degenerate group.  So, suppose that $E$ is a finite end of $M_\Gamma$ and that $\Phi$ is degenerate.   We can view $E$ as an end of $M_\Phi$ as well.   Let $E'$ be the complementary end of $E$ in $M_\Phi$.  Since $\Phi$ is degenerate, $E'$ is infinite, and so the end covering theorem implies that the restriction of the covering map  $\pi: M_\Phi\rightarrow M_\Gamma$ to $E'$ is finite-to-one.   Since the image of $E$ (viewed as an end of $M_\Phi$) under $\pi$ is just $E$ (viewed as an end of $M_\Gamma$),  the image of $E'$ in $M_\Gamma$ is the complementary end of $E$ in $M_\Gamma$.  Since the restriction of $\pi$ to $E'$ is finite-to-one, this implies that $\Phi$ has finite index in $\Gamma$, and so $\Gamma$ is isomorphic to the fundamental group of a surface, a contradiction.  

\medskip
\noindent
Since $E$ is either the infinite end of a degenerate group, or is a finite end of a quasifuchsian group, we see by Proposition \ref{bounded-notbounded-prop} that $E$ is bounded.
\end{proof}

\medskip
\noindent
In general, it is not possible to come up with a crisp statement of when
an end is bounded.  We present here a few examples to show what sorts of
things can go wrong.  There are some trivial situations in which an end is
always bounded.  For example, let $\Gamma$ be any non-elementary Kleinian group, and let $H$ be a small closed $3$-ball embedded in
$M_\Gamma$.  Then, under the definition we have given here, the complement
$E =M_\Gamma \setminus H$ is an end of $M_\Gamma$.  By definition, $E$ is
bounded as $H$ is compact.  

\medskip
\noindent
This example also illustrates why we have restricted our attention to peripheral ends.  Suppose that $M_\Gamma$ has several peripheral ends.  By taking the complement of a small compact set in $M_\Gamma$, we can construct an end of $M_\Gamma$ that contains several peripheral ends of $M_\Gamma$.   The behavior of such an end then becomes very complicated from the interaction of the several peripheral ends it contains.

\medskip
\noindent
Now, consider the following example of a peripheral end with
compressible boundary.  Let $\Phi_1$ be a purely loxodromic quasifuchsian
group, and let $\Phi_2$ be a purely loxodromic degenerate group.  Let
$\Gamma$ be the Klein combination of $\Phi_1$ and $\Phi_2$ (see Maskit
\cite{maskit}). Then, there is a component $\Delta$ of $\Omega(\Gamma)$
that is invariant under all of $\Gamma$.  The peripheral end of $M_\Gamma$
corresponding to $\Delta/\Gamma$ is then a finite end of $M_\Gamma$.  This end is not bounded, since its complementary piece contains the peripheral end
corresponding to the degenerate group.  Note that the end group
corresponding to this end is the whole group $\Gamma$.  However, using the
Klein-Maskit combination theorems, it is possible to realize this example
in more general Kleinian groups.

\medskip
\noindent
Finally, we can form a Kleinian group $\Gamma$ that is the free product of the
quasifuchsian group $\Phi_1$ and the degenerate group $\Phi_2$ in such a
way that $\Gamma$ contains no parabolic elements and so that every
component of $\Omega(\Gamma)$ is invariant under a conjugate of $\Phi_1$.  
(This Kleinian group is similar to the one constructed in the previous
paragraph, except that the invariant component has degenerated and is no
longer visible in the domain of discontinuity.  This construction is adapted from a construction due
to Maskit; see Section 5 of \cite{maskit-smash}.)  In this case,
$M_\Gamma$ has three peripheral ends: one finite end with incompressible
boundary, corresponding to the single surface $\Omega(\Gamma)$; one
infinite end with incompressible boundary, coming from the degenerate
group $\Phi_2$, and one infinite end with compressible boundary.  Note
that the end group corresponding to the infinite end with compressible
boundary is the whole group $\Gamma$, and this infinite end with
compressible boundary is not bounded.  As in the previous paragraph, it is
possible to realize this example in more general Kleinian groups using the
Klein-Maskit combination theorems.

\medskip
\noindent
We close this section by stating a mild refinement of the
tripartite division of points of the sphere at infinity ${\bf S}^2$ of
${\bf B}^3$ as given in Lemma \ref{tukia-lemma-1}.  The proof of Lemma \ref{tripartite} is essentially contained in the proofs of Lemma \ref{tukia-lemma-1} and Proposition \ref{incompressible-bounded}.  The main distinction
is that for topologically tame Kleinian groups, we are able to remove the
assumption of boundedness of the ends.

\begin{lemma}
Let $\Gamma$ be a purely loxodromic, topologically tame Kleinian group acting on ${\bf B}^3$.

\medskip
\noindent
If $\Gamma$ is doubly degenerate, let $E_0^0$ and $E_1^0$ be peripheral ends of $\Gamma$, and note that the associated end groups of $E_0^0$ and $E_1^0$ are both $\Gamma$.  Let $z\in {\bf S}^2$. Then either $z\in \Lambda_c(\Gamma)$, $z\in \Lambda_e(E_0^0)$, or $z\in \Lambda_e(E_1^0)$.

\medskip
\noindent
Suppose that $\Gamma$ is not doubly degenerate, and let $E^0$ be an incompressible infinite peripheral end of $\Gamma$.   Let $\Phi$ be the end group of $\Gamma$ associated to $E^0$, so that $\Phi$ is degenerate. Let $z\in {\bf S}^2$. Then either $z\in \Omega (\Phi)$,
$z\in \Lambda_c(\Phi)$ or $z$ is an end limit point of $E^0$. If
$z\in \Lambda (\Gamma)$ and $z$ is an endpoint of $E^0$, 
then $z\in \Lambda_e (E^0)$.
\label{tripartite}
\end{lemma}

\footnotesize{

}

\medskip
\noindent
{\footnotesize corresponding author:\\J. W. ANDERSON\\SCHOOL OF MATHEMATICS\\UNIVERSITY OF SOUTHAMPTON\\SOUTHAMPTON SO17 1BJ\\ENGLAND\\E-MAIL: j.w.anderson@maths.soton.ac.uk}

\medskip
\noindent
{\footnotesize K. FALK\\MATHEMATICAL INSTITUTE\\UNIVERSITY OF BERN\\SIDLERSTRASSE 5\\CH-3012 BERN\\SWITZERLAND\\E-MAIL: kurt.falk@math-stat.unibe.ch}

\medskip
\noindent
{\footnotesize P. TUKIA\\DEPARTMENT OF MATHEMATICS AND STATISTICS\\P. O. BOX 68 (GUSTAF H\"ALLSTR\"OMIN KATU 2B)\\FI-00014 UNIVERSITY OF HELSINKI\\FINLAND\\E-MAIL: pekka.tukia@helsinki.fi}


\begin{thebibliography}{99}

\bibitem{aaronson-sullivan} J. Aaronson and D. Sullivan, Rational
ergodicity of geodesic flows, {\em Erg. Th. Dyn. Syst.} {\bf 4}
(1984), 165--178.

\bibitem{agol} I. Agol, Tameness of hyperbolic 3-manifolds, preprint
(arXiv: math.GT/0405568).

\bibitem{ahlfors} L. V. Ahlfors, Finitely generated Kleinian groups, {\em
Amer. J. of Math.} {\bf 86} (1964), 413--29.

\bibitem{anderson} J. W. Anderson, The Limit Set Intersection Theorem for Finitely Generated Kleinian Groups, {\em Math. Research Lett.} {\bf 3} (1996), 675--692.

\bibitem{beardon} A. F. Beardon, The exponent of convergence of
Poincar\'e series, {\em Proc. London Math. Soc.} {\bf 18} (1968),
461--483.

\bibitem{beardon-book} A. F. Beardon, {\em The Geometry of Discrete Groups},
Springer Verlag, New York, 1983.

\bibitem{bijo} C. J. Bishop and P. W. Jones, The law of the iterated
logarithm for Kleinian groups, {\em Cont. Math.} {\bf 211} (1997), 17--50.

\bibitem{bonahon} F. Bonahon, Bouts des vari\'et\'es hyperboliques de
dimension 3, {\em Annals Math.} {\bf 124} (1986), 71--158.

\bibitem{calegari-gabai} D. C. Calegari and D. Gabai, 
Shrinkwrapping and the taming of hyperbolic $3$-manifolds, 
preprint (arXiv:math.GT/0407161)

\bibitem{canary-covering} R. D. Canary,
A Covering Theorem for Hyperbolic $3$-manifolds and its Applications,
{\em Topology} {\bf 35} (1996), 751--778.

\bibitem{canary2} R. D. Canary, Ends of hyperbolic $3$-manifolds, {\em J.
Amer. Math. Soc.} {\bf 6} (1993), 1--35.

\bibitem{canary1} R. D. Canary, Geometrically Tame Hyperbolic $3$-Manifolds,
{\em Proceedings of Symposia in Pure Mathematics} {\bf 54} (1993), 99--109.

\bibitem{culler-shalen} M. Culler and P. B. Shalen, Paradoxical
decompositions, $2$-generator Kleinian groups, and volumes of hyperbolic
$3$-manifolds, {\em J. Amer. Math. Soc.} {\bf 5} (1992), 231--288.

\bibitem{fatu} K. Falk and P. Tukia, A note on Patterson measures, 
preprint (2005).

\bibitem{hempel} J. Hempel, {\em $3$-manifolds}, {\em Annals of Mathematics Studies} {\bf 86}, Princeton University Press, 1976.

\bibitem{maskit} B. Maskit, {\em Kleinian groups}, Springer-Verlag, 1988.

\bibitem{maskit-smash} B. Maskit, Intersections of component 
subgroups of Kleinian Groups, in {\em Discontinuous Groups and 
Riemann Surfaces}, edited by L. Greenberg, 
{\em Annals of Mathematics Studies} {\bf 79}, 
Princeton University Press, 1974, 349--367.

\bibitem{morgan-bass} J.W.Morgan and H.Bass (editors), {\em The Smith
conjecture}, Academic Press 1984.

\bibitem{nicholls} P. J. Nicholls, {\em The Ergodic Theory of Discrete
Groups}, London Mathematical Society Lecture Notes Series, volume 143,
Cambridge University Press, 1989.

\bibitem{patterson} S. J. Patterson, The limit set of a Fuchsian
group, {\em Acta Math.} {\bf 136} (1976), 241--273.

\bibitem{patterson4} S. J. Patterson, Lectures on measures on limit sets of
Kleinian groups, {\em Analytical and geometric aspects of hyperbolic
space}, editor D.B.A. Epstein, Cambridge University Press, Cambridge, 1987.

\bibitem{rees1} M. Rees, Checking ergodicity of some geodesic flows with
infinite Gibbs measure, {\em Erg. Th. Dyn. Syst.} {\bf 1} (1981) 107--133.

\bibitem{rees2} M. Rees, Divergence type of some subgroups
of finitely generated Fuchsian groups, {\em Erg. Th. Dyn. Syst.} {\bf 1}
(1981) 209--221.

\bibitem{roblin} T. Roblin, Sur l'ergodicit\'e rationelle et les
propri\'et\'es ergodiques du flot g\'eod\'esique dans les vari\'et\'es
hyperboliques, {\em Ergod. Th. \& Dynam. Sys.} {\bf 20} (2000), 1785-1819.

\bibitem{scott-core} G. P. Scott,  Compact submanifolds of 3-manifolds, 
{\em J. London Math. Soc.} {\bf 7} (1973), 246--250.

\bibitem{sullivan1} D. Sullivan, The density at infinity of a discrete
group of hyperbolic motions, {\em Inst. Hautes \'Etudes Sci. Publ. Math.}
{\bf 50} (1979), 171-202.

\bibitem{sullivan2} D. Sullivan, Entropy, Hausdorff measures old and new,
and limit sets of geometrically finite Kleinian groups, {\em Acta. Math.}
{\bf 153} (1984), 259--277.

\bibitem{sullivan-growth} D. Sullivan, Growth of positive harmonic
functions and Kleinian group limit sets of zero planar measure and
Hausdorff dimension two, in Lecture Notes in Mathematics {\bf 894},
Springer Verlag, 1980, 127--144.

\bibitem{sullivan-positivity} D. Sullivan, Related aspects of positivity
in Riemannian geometry, {\em J. Difff. Geom.} {\bf 25} (1987), 327--351.

\bibitem{sullivan4} D. Sullivan, Travaux de Thurston sur les groupes
quasi-Fuchsiens et les varietes hyperboliques de dimension 3 fibrees
sur $S^1$, in {\em S\'eminaire Bourbaki, 32e ann\'ee, Vol. 1979/80, Exp.554},
number 842 in {\em Lect. Notes Math.}, Springer Verlag 1981, 196--214.

\bibitem{thurston} W.P. Thurston, {\em The geometry and topology of
three-manifolds}, Lecture notes from Princeton University 1978--1980.

\bibitem{tukia2} P. Tukia, The Hausdorff dimension of the limit set of a
geometrically finite Kleinian group, {\em Acta. Math.} {\bf 152} (1984),
127--140.

\bibitem{tukia3} P. Tukia, On isomorphisms of geometrically finite
M\"obius groups, {\em Publ.  Math. IHES} {\bf 61} (1985), 171--214.

\bibitem{tukia1} P. Tukia, The Poincar\'e series and the conformal measure
of conical and Myrberg limit points, {\em J. Anal. Math.} {\bf 62} (1994),
241--259.

\end{thebibliography}
\end{document}